\documentclass[11pt,letter]{amsart}

\usepackage{amssymb}
\usepackage{color}
\usepackage{graphicx}
\usepackage{float}
\usepackage[all,cmtip]{xy}
\usepackage{tikz}
\usetikzlibrary{matrix}
\usepackage{url}
\usepackage{subfig}
\usepackage{hyperref}
\hypersetup{hidelinks}
\usepackage{comment}
\usepackage{appendix}
\mathchardef\mhyphen="2D

\newcommand{\std}{{\operatorname{std}}}

\newcommand{\Diff}{\operatorname{Diff}}
\newcommand{\Ort}{\operatorname{O}}
\newcommand{\SO}{\operatorname{SO}}
\newcommand{\U}{\operatorname{U}}

\newcommand{\SU}{\operatorname{SU}}

\newcommand{\rel}{\operatorname{rel}}

\newcommand{\even}{\operatorname{even}}

\newcommand{\Id}{{\operatorname{Id}}}

\newcommand{\rot}{\operatorname{rot}}
\newcommand{\tb}{\operatorname{tb}}

\newcommand{\Cont}{\operatorname{Cont}}

\newcommand{\Fr}{\operatorname{Fr}}
\newcommand{\CFr}{\operatorname{CFr}}
\newcommand{\ev}{\operatorname{ev}}

\newcommand{\R}{{\mathbb{R}}}
\newcommand{\C}{{\mathbb{C}}}
\newcommand{\Z}{{\mathbb{Z}}}

\newcommand{\NS}{{\mathbb{S}}}

\newcommand{\D}{{\mathbb{D}}}

\newcommand{\PP}{{\mathbb{P}}}
\newcommand{\Op}{{\mathcal{O}p}}

\newcommand{\Lag}{{\mathfrak{Lag}}}
\newcommand{\Leg}{{\mathfrak{Leg}}}

\newcommand{\Emb}{{\mathfrak{Emb}}}

\newcommand{\CEmb}{{\mathfrak{CEmb}}}
\newcommand{\CStr}{{\mathfrak{C}Str}}
\newcommand{\Trans}{{\mathfrak{Trans}}}

\newtheorem{theorem}{Theorem}[subsection]
\newtheorem{lemma}[theorem]{Lemma}
\newtheorem{proposition}[theorem]{Proposition}
\newtheorem{corollary}[theorem]{Corollary}

\theoremstyle{definition}
\newtheorem{definition}[theorem]{Definition}

\newtheorem{remark}[theorem]{Remark}

\makeatletter
\@namedef{subjclassname@2020}{\textup{2020} Mathematics Subject Classification}
\makeatother

\begin{document} 

\title[Spaces of tight contact structures in dimension 3]{Spaces of tight contact structures in dimension 3}

\subjclass[2020]{Primary: 53D35, 57K33.}

\keywords{}

\author[E.Fern\'{a}ndez]{Eduardo Fern\'{a}ndez}
\address{School of Mathematics\\ Georgia Institute of Technology\\ Atlanta\\ GA, USA}
\email{eduardofernandez@gatech.edu}

\author[J.Mart\'{i}nez-Aguinaga]{Javier Mart\'{i}nez-Aguinaga}
\address{Universidad Complutense de Madrid, Departamento de Álgebra, Geometría y Topología,
	Facultad de Ciencias Matem\'{a}ticas. Plaza de las Ciencias, 3, 28040 Madrid, Spain.}
\email{frmart02@ucm.es}

\author[F.Presas]{Francisco Presas}
\address{Instituto de Ciencias Matem\'{a}ticas CSIC-UAM-UC3M-UCM, C. Nicol\'{a}s Cabrera, 13-15, 28049 Madrid, Spain.}
\email{fpresas@icmat.es}

\begin{abstract}
    We study spaces of tight contact structures and contactomorphism groups in dimension $3$. We introduce a microfibration method for studying spaces of convex surfaces in contact $3$-manifolds and apply it to describe spaces of convex disks and spheres in tight contact $3$-manifolds. With these ingredients in hand, we prove a parametric gluing theorem for contact handle attachments: the weak homotopy type of the space of tight contact structures is invariant under contact $0$-, $1$-, and $3$-handle attachments. 

    Based on these new techniques, we present a substantial number of independent applications. First, we obtain explicit descriptions of the weak homotopy type of several spaces of tight contact structures. We determine the weak homotopy type of the contactomorphism group of every standard contact handlebody, of the standard $\NS^1\times\NS^2$, and of every Legendrian circle bundle with non-empty boundary, partially resolving a conjecture of Giroux. Additionally, we compute the weak homotopy type of an infinite family of spaces of Legendrian and transverse unknots, proving Legendrian and transverse analogues of Hatcher's theorem for smooth unknots in the $3$-sphere.
\end{abstract}

\maketitle
\setcounter{tocdepth}{1} 
\let\thefootnote\relax\footnotetext{The authors acknowledge the use of OpenAI's ChatGPT solely for proofreading, identifying errors and typos, and improving the English of the manuscript. AI was not used to generate any mathematical content, arguments, or results.}


\section{Introduction}\label{sec:intro}


\subsection{Context}

While classical contact topology is primarily concerned with the classification of contact structures, this article initiates a systematic study of parametrized contact topology in dimension $3$, namely, the study of the weak homotopy type of spaces of contact structures. Although the present article is concerned exclusively with tight contact structures \cite{EliashbergOT}, many of the techniques developed here extend well beyond this setting and have recently been applied to the study of spaces of overtwisted contact structures \cite{F:SOT}.

Conceptually, our approach to the study of spaces of tight contact structures is simple: the classification of tight contact structures in dimension three is largely based on cutting contact manifolds along convex surfaces and understanding the resulting pieces; see, for instance, \cite{CGH,Colin,GirouxInventiones,Giroux,GirouxCriterion,Ka1997,Honda,honda2000b,HondaGluing,HondaKazezMatic}. This program relies heavily on Giroux's convex surface theory \cite{GirouxConvex}, Honda's theory of bypass attachments \cite{Honda}, gluing results for tight contact structures along convex surfaces \cite{Colin,HondaGluing}, and the uniqueness of tight contact structures on the $3$-disk \cite{Eliashberg}. 

Since Giroux's convex surface theory has not been developed in a genuinely parametric setting beyond $0$- and $1$-parameter families of surfaces, results concerning contactomorphism groups and spaces of tight contact structures remain comparatively scarce, highly case-dependent, and usually do not extend beyond computations of the fundamental group \cite{Bourgeois,DingGeiges,GeigesGonzalo04,GeigesKlukas,Giroux,GirouxMassot,HyunkiMin,MunozEchaniz}. One notable exception is the theorem of Eliashberg and Mishachev \cite{EliashbergMishachevTight}, which asserts that the space of tight contact structures on the $3$-disk is contractible. This result serves as a building block throughout the article.

Our goal is to lift this cut-and-paste philosophy to a genuinely parametric setting, thereby making it applicable to the study of higher homotopy groups. We achieve this by introducing a new technique for studying certain higher-dimensional families of convex surfaces, which we call the \emph{embedding microfibration method} (Theorem \ref{thm:naked}).

This method leads to the first complete descriptions of the weak homotopy type of various spaces of convex surfaces and, as a consequence, to the first gluing results for families of tight contact structures. It also yields several complete computations of the weak homotopy type of spaces of tight contact structures and contactomorphism groups of tight contact $3$-manifolds, as well as the first complete descriptions of various Legendrian and transverse embedding spaces in tight contact $3$-manifolds. Since the first version of this article appeared, the method has found further applications, including rigidity phenomena in contact mapping class groups, new parametric gluing results, structural properties of Legendrian embedding spaces under topological knot operations, and the study of spaces of overtwisted contact structures \cite{FME22,FernandezMin,F:SOT}.

\subsection{Spaces of convex surfaces} 

To ease the reading, we defer the statement of our microfibration method (Theorem \ref{thm:naked}) to later in the article and proceed directly to the applications that motivate it. The main outcome is a description of the weak homotopy type of several spaces of convex surfaces. 

\subsubsection{Convex disks with Legendrian boundary}

Let $(M,\xi)$ be a tight contact $3$-manifold. Let $e:\D^2\hookrightarrow(M,\xi)$ be a convex embedding whose boundary $\partial\D^2$ is a Legendrian unknot satisfying
\[
\tb(\partial\D^2)+|\rot(\partial\D^2)|=-1.
\]
Let $\mathcal{F}:=e^*\xi$ denote the characteristic foliation induced on $\D^2$ by the embedding. See Figure \ref{fig:standard} for an example of such a characteristic foliation in the case $\tb=-1$  (a convex disk carrying this characteristic foliation is referred to as \emph{standard}).

Let $\Emb(\D^2,M)$ denote the space of embeddings of disks that agree with $e$ on an open neighborhood of the boundary. We denote by $\Emb((\D^2,\mathcal{F}),(M,\xi))$ the subspace of $\Emb(\D^2,M)$ consisting of convex embeddings $j:\D^2\hookrightarrow(M,\xi)$ that induce the same characteristic foliation as $e$, namely $j^*\xi=\mathcal{F}$. Our first result is the following.

\begin{theorem}\label{thm:DisksIntro}
Let $(M,\xi)$ be a tight contact $3$-manifold. Then the natural inclusion
\[
\Emb((\D^2,\mathcal{F}),(M,\xi))\hookrightarrow \Emb(\D^2,M)
\]
is a weak homotopy equivalence.
\end{theorem}

\begin{remark}
The analogous statement also holds for embeddings of disks bounding a positively transverse unknot with self-linking number $-1$; see Theorem \ref{thm:MiniDisks}.
\end{remark}

\begin{remark}
The $\pi_0$-surjectivity of the inclusion follows from the work of Giroux \cite{GirouxConvex,GirouxCriterion}, while $\pi_0$-injectivity follows from the work of Colin \cite{Colin}. See also \cite{Eliashberg,EliashbergTransverse,EliashbergFraser}.
\end{remark}

\subsubsection{Standard convex spheres}

Let $e:\NS^2\hookrightarrow (M,\xi)$ be a convex embedding of a sphere into a tight contact $3$-manifold, and let $\mathcal{F}=e^*\xi$ denote its characteristic foliation. An important example is the characteristic foliation induced on the boundary of a Darboux ball, which has exactly two elliptic singularities, one at each pole. Convex spheres carrying this precise characteristic foliation are referred to as \emph{standard}. Let $\Emb(\NS^2,M)$ denote the space of smooth embeddings of spheres into $M$, and let $\Emb((\NS^2,\mathcal{F}),(M,\xi))$ denote the subspace consisting of embeddings that induce the same characteristic foliation as $e$.

\begin{theorem}\label{thm:SpheresIntro}
Let $(M,\xi)$ be a tight contact $3$-manifold. Then, for every $k>0$, there is an isomorphism
\[
\pi_k\big(\Emb(\NS^2,M),\Emb((\NS^2,\mathcal{F}),(M,\xi))\big)
\cong
\pi_k(\SO(3),\U(1)).
\]
\end{theorem}

The isomorphism is induced by evaluating the $1$-jet map at a point. This result may be viewed as a natural parametric analogue of Colin's theorem \cite{Colin}.

\subsection{Gluing Theorems: Contact Handle Attachments} Giroux proved that every contact manifold admits a contact handle decomposition \cite{GirouxICM}, providing a positive answer to a question posed by Eliashberg and Gromov \cite{EliashbergGromov:Convex}. 

Informally, for contact $3$-manifolds, the description of contact handles is simple: the $0$- and $3$-handles are Darboux balls, while the $1$- and $2$-handles are $I$-invariant convex neighborhoods of standard disks; see \cite{HH19,Ozbagci} for the precise definitions. Later in the manuscript, we will define contact handle attaching maps at the level of spaces of contact structures; see Section \ref{sec:Gluing}.

It is clear that contact $0$- and $3$-handle attachments preserve tightness. Moreover, a celebrated theorem of Colin \cite{Colin} states that contact $1$-handle attachments also preserve tightness. Together, \cite{Colin} and \cite{Eliashberg} imply that contact $0$-, $1$-, and $3$-handle attachments induce a one-to-one correspondence between isotopy classes of tight contact structures (rel. boundary). On the other hand, contact $2$-handle attachments do not preserve tightness in general, as not every contact $3$-manifold is tight.

The importance of Theorems \ref{thm:DisksIntro} and \ref{thm:SpheresIntro} is that they allow us to prove parametric gluing results for tight contact structures. The precise formulation of our gluing theorem for contact handle attachments will be given later in Theorem \ref{thm:Gluing}. For the purposes of this introduction, we state the following non-technical version, which captures the main geometric content of the result:

\begin{theorem}\label{thm:ContactHandleIntro}
The weak homotopy type of the space of tight contact structures is invariant under contact $0$-, $1$-, and $3$-handle attachments. 
\end{theorem}

\begin{remark}
    For simplicity, we suppress a minor technical subtlety in the  $3$-handle attachment case. See Theorem \ref{thm:Gluing} for the precise statement. 
\end{remark}

\begin{remark}
The fact that the contact $0$-handle attachment map is a weak homotopy equivalence follows from Eliashberg--Mishachev \cite{EliashbergMishachevTight}; i.e., the space of tight contact structures on the $3$-ball is contractible. In this article, we establish the corresponding statement for the contact $1$- and $3$-handle attachments.
\end{remark}

In particular, any nontrivial change in the homotopy type of the space of tight contact structures must arise from contact $2$-handle attachments. Assuming that tightness is preserved, the main obstruction to extending the previous statement to $2$-handles is that the attaching region is an annulus, and we are unable to establish general theorems for such spaces of convex annuli. Such annuli are called \emph{standard} in the work of Min and the first author \cite{FernandezMin}. In the present article, we study these annuli in the special case of Legendrian circle bundles. 

\begin{remark}Using the methods developed here, several gluing results along annuli, analogous to Theorem \ref{thm:ContactHandleIntro}, were established in \cite{FernandezMin}. Moreover, gluing results for connected sums were established in \cite{FME22} using Theorem \ref{thm:SpheresIntro}, yielding a parametric version of Colin's Decomposition Theorem \cite{Colin}.
\end{remark}

\subsection{Spaces of tight contact structures and contactomorphism groups}

We now turn to specific cases in which the previous results allow us to completely determine the weak homotopy type of the space of tight contact structures. Through Gray stability \cite{GrayStability}, these computations also yield information about contactomorphism groups. 

Indeed, for every compact contact $3$-manifold $(M,\xi)$ there is a fibration
\[
\Cont(M,\xi)\cap \Diff_0(M)\hookrightarrow \Diff_0(M)\longrightarrow \CStr(M,\xi),
\]
which we refer to as the \emph{Gray fibration}; see \cite{GirouxMassot}. Here, $\Diff_0(M)$ denotes the path-connected component of the identity in $\Diff(M)$, the group of diffeomorphisms of $M$ that are the identity near the boundary, $\Cont(M,\xi)\subseteq \Diff(M)$ is the subgroup consisting of those diffeomorphisms that preserve the contact structure $\xi$, known as \emph{contactomorphisms}, and $\CStr(M,\xi)$ denotes the path-connected component of $\xi$ in the space of contact structures $\CStr(M)$ that coincide with $\xi$ near the boundary. 

The homotopy type of diffeomorphism groups of many $3$-manifolds is relatively well understood \cite{bamler,bamler2,HatcherSmale,HatcherSpheres,HatcherHomeo,EllipticDiff,HendriksLaudenbach,Gabai}. Consequently, we will regard the homotopy type of $\Diff(M)$ as known input, and the homotopy groups of spaces of contact structures are closely related to understanding those of contactomorphism groups. Throughout the article, we will freely move between these two points of view.

\subsubsection{Standard contact handlebodies}

Let $(H_g,\xi)$ be a genus $g$-handlebody equipped with a tight contact structure. We say that $(H_g,\xi)$ is \emph{standard} if it admits $g$ pairwise disjoint non-separating disks whose boundaries are Legendrian unknots with $\tb=-1$. Equivalently, standard contact handlebodies are obtained from a Darboux ball by contact $1$-handle attachments. Our gluing result Theorem \ref{thm:ContactHandleIntro} can then be applied to show that:

\begin{theorem}\label{thm:StandardHandlebodyINTRO}
Let $(H_g,\xi)$ be a standard contact handlebody. The space $\CStr(H_g,\xi)$ is weakly contractible. Moreover, the inclusion $\Cont(H_g,\xi)\hookrightarrow \Diff(H_g)$ is a weak homotopy equivalence. In particular, $\Cont(H_g,\xi)$ is weakly contractible.
\end{theorem}

\begin{remark} As a consequence, the space of tight contact structures on $\NS^1\times\R^2$ which coincide with $(J^1 \NS^1,\xi_\std)$ at infinity is weakly contractible (see Corollary \ref{cor:J1S1}). 
\end{remark}

\begin{remark}
    In Theorem \ref{thm:MultistandardHandlebody} we will prove an analogous statement for a wider class of (universally) tight handlebodies, which we call \emph{multistandard}. 
\end{remark}

\subsubsection{The standard contact $\NS^1\times\NS^2$}

Let $(\NS^1\times\NS^2,\xi_\std)$ be the standard tight contact structure on $\NS^1\times \NS^2$, that is, the $\NS^1$-invariant contact structure over a standard convex sphere $\NS^2$. This is the unique, up to isotopy, tight contact structure on $\NS^1\times \NS^2$ \cite{Eliashberg}. In coordinates, $\NS^1\times\NS^2\subseteq \NS^1\times \R^3(\theta,(x,y,z))$ we have that $\xi_\std=\ker (xd\theta+ydz-zdy)$. 

To describe the space $\CStr(\NS^1\times\NS^2,\xi_\std)$, of tight contact structures on $\NS^1\times\NS^2$, we introduce the following evaluation map. Fix a global framing of $\NS^1\times\NS^2$, so that it coincides with $\langle \partial_\theta, \partial_y, \partial_z \rangle $ along $\NS^1\times \{(1,0,0)\}$. 

Since every contact structure $\xi\in \CStr(\NS^1\times\NS^2,\xi_\std)$ is cooriented, it has a well-defined positively oriented unit normal vector field $R^\xi$. We thus obtain a well-defined \emph{evaluation map} along $\NS^1\equiv \NS^1\times\{(1,0,0)\}$ by restriction of this vector field along $\NS^1$: \begin{equation}\label{eq:EvaluationMapS1xS2}
\ev_{\NS^1}:\CStr(\NS^1\times\NS^2,\xi_\std) \rightarrow L\NS^2, \qquad \xi\mapsto R^\xi_{|\NS^1\times \{(1,0,0)\}}
\end{equation}

\begin{theorem}\label{thm:ContactomorphismsS1xS2} The following statements hold.
\begin{itemize}
\item [(i)] The evaluation map (\ref{eq:EvaluationMapS1xS2}) is a weak homotopy equivalence.
\item [(ii)] The group $\Cont(\NS^1\times\NS^2,\xi_\std)$ is weakly homotopy equivalent to $\NS^1\times\Ort(2)\times\Omega \U(1)$.
\end{itemize}
\end{theorem}

\subsubsection{Legendrian circle bundles over compact orientable surfaces with non-empty boundary}

The study of the path-connected components of the contactomorphism group of a Legendrian circle bundle over a compact orientable surface was initiated by E. Giroux \cite{Giroux} and completed by E. Giroux and P. Massot \cite{GirouxMassot}. A Legendrian circle bundle over a surface is a contact $3$-manifold $(V,\xi)$ equipped with an $\NS^1$-bundle structure over a surface $S$ and such that $\xi$ is tangent to the fibers. The paradigm to keep in mind is the space of cooriented contact elements $(\NS(T^*S),\xi_\std)$.

As observed by Lutz, every Legendrian circle bundle is a fibered contact covering of a space of cooriented contact elements \cite{Lutz}. Thus, there is an inclusion of the diffeomorphism group of the base $\Diff(S)$ into the contactomorphism group $\Cont(V,\xi)$ of $(V,\xi)$. In \cite[Corollary 2.8]{GirouxMassot} the authors prove that this inclusion is a $\pi_0$-isomorphism if the boundary of $S$ is non-empty, see also \cite{Giroux}. In Giroux's article \cite{Giroux}, it is written: 
\begin{quote}
   \em ``En fait, les plongements ci-dessus sont tr\`es probablement des \'equivalences d’homotopie. Au prix de complications surtout techniques, la m\'ethode suivie ci-apr\`es semble d’ailleurs applicable aux familles – de diff\'eomorphismes, de plongements, de structures de contact – d\'ependant d’un nombre quelconque de param\`etres.'' \em
\end{quote}
This conjecture of Giroux was one of the starting points of this article. In fact, we are able to prove the conjecture for the case when the base has non-empty boundary:

\begin{theorem}\label{thm:LegendrianBoundary}
Let $(V,\xi)$ be a Legendrian circle bundle over a compact orientable surface $S$ with non-empty boundary. The natural inclusion \begin{equation}\label{eq:InclusionCircleBundles} i:\Diff(S)\hookrightarrow\Cont(V,\xi)\end{equation} is a weak homotopy equivalence. Moreover, $\CStr(V,\xi)$ is a weakly contractible space.
\end{theorem}

\subsection{Legendrian and transverse embedding spaces} The new techniques developed in this article allow us to tackle open questions in the study of Legendrian and transverse embedding spaces and, in particular, their fundamental groups  \cite{FMP,Kalman}. 

We will denote by $\Leg(M,\xi)$ the space of Legendrian embeddings into $(M,\xi)$. For a Legendrian embedding $\Lambda:\NS^1\hookrightarrow (M,\xi)$, we will denote by $\Leg(M,\xi,\Lambda)\subseteq \Leg(M,\xi)$ the path-connected component containing $\Lambda$, and by \[(C(\Lambda,M),\xi):=(M\backslash \Op(\Lambda),\xi)\] its knot exterior. The contact isotopy extension theorem implies the existence of a fibration
\[ \Cont(C(\Lambda,M),\xi) \hookrightarrow \Cont(M,\xi)\rightarrow \Leg (M,\xi), \]
induced by postcomposition. Therefore, the homotopical study of Legendrian embedding spaces reduces to the study of contactomorphism groups and, as discussed above, the latter to the study of spaces of contact structures. For transverse embedding spaces, there are similar fibration sequences.

\begin{remark}
    This strategy to study Legendrian embedding spaces mimics the strategy of Hatcher, further refined by Budney, to study the homotopy type of knot spaces in the $3$-sphere \cite{Budney,HatcherSmale,HatcherKnots1,HatcherKnots2}. Their approach is based on the Cerf-Palais fibration theorem \cite{Cerf:Fibration,Palais}, i.e. the parametrized isotopy extension theorem.
\end{remark}

\begin{remark}
    In fact, the relation between contactomorphism groups and Legendrian embedding spaces is even deeper. The entire weak homotopy type of the contactomorphism group can be recovered from the space of Legendrian embeddings of a suitable Legendrian \emph{graph} (see Corollary \ref{cor:GirouxSkeleton}).
\end{remark}

\subsubsection{Legendrian and transverse unknots} 
\begin{figure}[!h]
  \centering
{\includegraphics[width=0.65\textwidth]{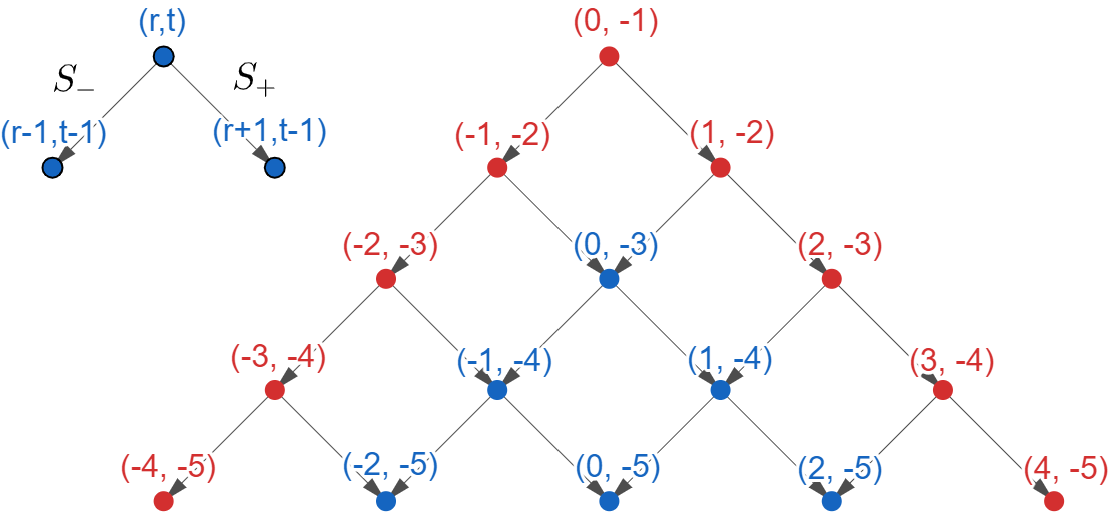}}
  \caption{Legendrian mountain range depicting the path-connected components of the space of Legendrian unknots. The boundary components to which Theorem \ref{thm:LegendrianAndTransverseUnknotsINTRO} applies are shown in red; the remaining components are shown in blue. \label{fig:LegendrianUnknots}}
\end{figure}

We begin our study of Legendrian and transverse embedding spaces by considering unknots in tight contact $3$-manifolds. The path-connected components of the space of oriented Legendrian unknots in a tight contact $3$-manifold are well understood, due to the Eliashberg--Fraser Theorem \cite{EliashbergFraser}. Every Legendrian unknot is obtained from the unique Legendrian unknot with $\tb=-1$ by a finite sequence of positive and negative stabilizations. See Figure \ref{fig:LegendrianUnknots}. Similarly, transverse unknots are classified by their self-linking number and are all obtained from the unique transverse unknot with self-linking number $-1$, see \cite{EliashbergTransverse}.

In order to state our main result, we introduce some notation. Let $(M,\xi)$ be a tight contact $3$-manifold, and fix a unit vector $(p,v)\in \NS(T_pM)$. For a smooth unknot $U:\NS^1\hookrightarrow M$, we denote by $\Emb_{(p,v)}(\NS^1,M,U)$ the space of \emph{long unknotted embeddings}. Here, \emph{long} means that every unknotted  embedding $\gamma:\NS^1\hookrightarrow M$ in this space satisfies $\gamma(0)=p$ and $\gamma'(0)=v$.

If $\Lambda_{t,r}:\NS^1\hookrightarrow (M,\xi)$ is a Legendrian unknot with $\tb(\Lambda_{t,r})=t$ and $\rot(\Lambda_{t,r})=r$, and $(p,v)\in \NS(\xi_p)\subseteq \NS(T_pM)$ is a Legendrian vector, we also denote by \[ \Leg_{(p,v)}(M,\xi,\Lambda_{t,r}):= \Leg(M,\xi,\Lambda_{t,r})\cap \Emb_{(p,v)}(\NS^1,M,\Lambda_{t,r}) \subseteq \Emb_{(p,v)}(\NS^1,M,\Lambda_{t,r}) \]  the subspace of \emph{long Legendrian unknots} with Thurston-Bennequin invariant $\tb=t$ and rotation number $\rot=r$.

Similarly, let $T_s:\NS^1\hookrightarrow (M,\xi)$ be a transverse unknot with self-linking number $\mathrm{sl}(T_s)=s$ and $(p,v)\in \NS(T_p M)$ a positively transverse vector. Then, $\Trans(M,\xi,T_s)$ denotes the path-connected component of the space of transverse embeddings that contains $T_s$, while 
\[ \Trans_{(p,v)}(M,\xi,T_s):= \Trans(M,\xi,T_s)\cap \Emb_{(p,v)}(\NS^1,M,T_s) \subseteq \Emb_{(p,v)}(\NS^1,M,T_s) \]  is the subspace of \emph{long transverse unknots} with self-linking  invariant $\mathrm{sl}=s$.

\begin{theorem}\label{thm:LegendrianAndTransverseUnknotsINTRO}
    Let $(M,\xi)$ be a tight contact $3$-manifold. 
    \begin{itemize}
        \item [(i)] For a Legendrian unknot $\Lambda_{t,r}$ such that $t+|r|=-1$, the natural inclusion 
    \[ \Leg_{(p,v)}(M,\xi,\Lambda_{t,r})\hookrightarrow \Emb_{(p,v)}(\NS^1,M,\Lambda_{t,r}) \] is a weak homotopy equivalence. 
       \item [(ii)] For the transverse unknot $T_{-1}$ with self-linking number $-1$, the natural inclusion     \[ \Trans_{(p,v)}(M,\xi,T_{-1})\hookrightarrow \Emb_{(p,v)}(\NS^1,M,T_{-1}) \] is a weak homotopy equivalence. 
       \end{itemize}
\end{theorem}

For Legendrian unknots the result will be deduced directly from Theorem \ref{thm:DisksIntro} by relating the space of Legendrian unknots with the space of convex disks with fixed Legendrian boundary; see Lemma \ref{lem:UnknotsVSDisks}. A similar argument, relating the study of smooth unknots to that of smooth disks via half-disk embeddings, appeared independently in the work of Budney and Gabai \cite{BudneyGabai}. For the transverse unknot a similar argument works. 

\begin{remark}
Every Legendrian knot $\Lambda$ admits a transverse push-off $T(\Lambda)$, called the \emph{transversalization} of $\Lambda$, which is unique up to transverse isotopy. Moreover,
\[ \mathrm{sl}(T(\Lambda))=\tb(\Lambda)-\rot(\Lambda).\]
Observe that the right-hand side is invariant under negative stabilization. In fact, it is well known that the classification of transverse knots is equivalent to the classification of Legendrian knots up to negative stabilization; see, for instance, \cite{Etnyre}. Although we do not pursue this perspective in the present article, a parametric version of this correspondence can be used to deduce part (ii) from part (i) of Theorem \ref{thm:LegendrianAndTransverseUnknotsINTRO}.
\end{remark}

 \begin{remark} Theorem \ref{thm:LegendrianAndTransverseUnknotsINTRO} is purely $3$-dimensional, as the work in preparation of Eliashberg and Kragh shows \cite{EliashbergKragh}: in the higher-dimensional standard $(\R^{2n+1},\xi_\std)$ there are exotic families of long Legendrian unknots. Here, exotic refers to formally contractible but not geometrically contractible, see \cite{MurphyLoose,FMP}.\end{remark}

\subsubsection{Legendrian and transverse unknots in the standard contact $3$-sphere}
Theorem \ref{thm:LegendrianAndTransverseUnknotsINTRO} almost completely determines the weak homotopy type of the space $\Leg(M,\xi,\Lambda_{t,r})$. Indeed, there exists a natural fibration map \[ \Leg_{(p,v)}(M,\xi,\Lambda_{t,r})\hookrightarrow \Leg(M,\xi,\Lambda_{t,r})\rightarrow \NS(\xi). \]
There is a similar fibration sequence for transverse knots. 

In the case of $(\NS^3,\xi_\std)$, this description is particularly simple. A \em Legendrian great circle \em is the intersection of a Lagrangian plane in $\C^2$ with $\NS^3$, whereas a \em transverse great circle is \em the intersection of a complex line in $\C^2$ with $\NS^3$. Observe that the space of parametrized Legendrian great circles is naturally identified with the Stiefel manifold of Lagrangian frames $V_{4,2}^{\Lag}\cong \U(2)\subseteq V_{4,2}$, and the space of parametrized transverse great circles with the Stiefel manifold of complex frames $V_{4,2}^{\C}\cong \SU(2)\subseteq V_{4,2}$.

The following result provides an affirmative answer to a question of Chantraine, Dimitroglou Rizell, Ghiggini, and Golovko \cite[Remark 4.2]{CDGG}. 

\begin{theorem}\label{thm:UnknotsS3Intro} The following statements hold:
    \begin{itemize} 
        \item [(i)] The natural inclusion \[ V_{4,2}^{\Lag}\hookrightarrow \Leg(\NS^3,\xi_\std,\Lambda_{-1,0}) \] is a weak homotopy equivalence. 
        \item [(ii)] The natural inclusion \[ V_{4,2}^{\C}\hookrightarrow \Trans(\NS^3,\xi_\std,T_{-1}) \] is a weak homotopy equivalence. 
    \end{itemize}
\end{theorem}

\begin{remark}
    This theorem is the contact analogue of a celebrated theorem of A. Hatcher in the smooth category, which asserts that the space of parametrized smooth unknots in $\NS^3$ is weakly homotopy equivalent to the space $V_{4,2}$ of parametrized great circles \cite{HatcherSmale,HatcherKnots1}; see also \cite{Budney}. Our proof makes essential use of Hatcher's theorem.
\end{remark}
\begin{remark}
Theorem \ref{thm:UnknotsS3Intro}(i), together with Eliashberg--Polterovich's result on the space of Lagrangian fillings of the standard Legendrian unknot \cite{EliashbergPolterovich}, provides a complete description of the weak homotopy type of the space of Lagrangian fillings of \emph{all} Legendrian unknots with $\tb=-1$, i.e. without fixing the boundary (see Corollary \ref{cor:Lagrangians}).
\end{remark}

\subsubsection{Legendrian $(n,n)$-torus links with max-$\tb$ number}

We conclude our investigation of Legendrian embedding spaces by studying the space  $\mathcal{L}_n$ of Legendrian embeddings of the maximal-$\tb$ invariant $(n,n)$-torus link into $(\NS^3,\xi_\std)$. 

\begin{theorem}\label{thm:LegendrianNNlinks}
    Let $\mathcal{M}_n$ denote the mapping class group of a $2$-sphere with $n$-holes. Then, there is a weak homotopy equivalence \[ \mathcal{L}_n\cong\U(2)\times K(\mathcal{M}_n,1).\]
\end{theorem}

\begin{remark} This result should be compared with the work of Casals and Gao \cite{CasalsGao} on the $(4,4)$-torus link. In that article, the authors use invariants from microlocal sheaf theory to construct a surjective homomorphism from the fundamental group of the space of \em unparametrized \em $(4,4)$-torus links onto the mapping class group of the $2$-sphere with $4$ marked points. 
\end{remark}

\begin{remark} Using the methods introduced in the present article, this result was subsequently extended to all max-$\tb$ positive torus links and, more generally, to a broader class of Legendrian links, including every max-$\tb$ representative of an algebraic link, in the work of Min and the first author \cite{FernandezMin}.
\end{remark}

\subsection{Outline}

The article is organized as follows. Sections \ref{sec:Basics}, \ref{sec:Hatcher}, and \ref{sec:Eliash} cover the preliminary material needed throughout the article. In Section \ref{sec:Basics}, we review several results from convex surface theory that will be used throughout the manuscript, as well as fibration sequences relating the main spaces considered here, including spaces of convex surfaces, contact structures, and Legendrian embeddings. Section \ref{sec:Hatcher} provides a brief review of some of Hatcher's results on diffeomorphism groups and embedding spaces in dimension three. This also serves as a useful warm-up, since several of the arguments discussed there will later be translated into the contact setting using the fibrations introduced in Section \ref{sec:Basics}. Section \ref{sec:Eliash} recalls the foundational work of Eliashberg--Mishachev \cite{EliashbergMishachevTight} on tight contact structures on the $3$-disk and $3$-sphere.

The remaining sections contain the main new contributions of the article and are organized as follows. In Section \ref{sec:MicrofibrationTrick}, we introduce and prove our \emph{embedding microfibration method} (Theorem \ref{thm:naked}). This is one of the main new techniques introduced in the article and will allow us to study several spaces of convex surfaces. In Section \ref{sec:ConvexSurfaces}, we give the first applications of the microfibration method. In particular, we prove Theorems \ref{thm:DisksIntro} and \ref{thm:SpheresIntro}, as well as an analogous result for disks with transverse boundary (Theorem \ref{thm:MiniDisks}).

In Section \ref{sec:Gluing}, we turn to our parametric gluing theorem for contact handle attachments. We state and prove the technical version of Theorem \ref{thm:ContactHandleIntro} in Theorem \ref{thm:Gluing}. We then deduce several consequences, including our result on standard contact handlebodies (Theorem \ref{thm:StandardHandlebodyINTRO}) and Theorem \ref{thm:UnknotsS3Intro}(i) on the standard Legendrian unknot in the $3$-sphere. We conclude the section with several structural consequences concerning contactomorphism groups and Giroux contact cell decompositions (Lemma \ref{lem:GirouxDecomposition} and Corollary \ref{cor:GirouxSkeleton}).

In Section \ref{sec:LegendrianUnknots}, we turn to the study of Legendrian unknots in a general tight contact $3$-manifold. We first prove Theorem \ref{thm:LegendrianAndTransverseUnknotsINTRO}(i). The main ingredient is Lemma \ref{lem:UnknotsVSDisks}, which identifies the space of convex Seifert disks with the based loop space of the corresponding space of Legendrian unknots. Combined with our results on convex disks, this yields the theorem. We conclude the section by indicating the analogous argument for transverse unknots with self-linking number $-1$, proving Theorems \ref{thm:LegendrianAndTransverseUnknotsINTRO}(ii) and \ref{thm:UnknotsS3Intro}(ii).

In Section \ref{sec:S1xS2}, we use our results on convex spheres to prove Theorem \ref{thm:ContactomorphismsS1xS2}. In Section \ref{sec:LegendrianFibrations}, we prove Theorem \ref{thm:LegendrianBoundary} by applying the microfibration method to certain spaces of fibered annuli. Theorem \ref{thm:LegendrianNNlinks} is proved in Section \ref{sec:TorusLinks}. Finally, in Section \ref{sec:TightCharts}, we prove an existence result for common trivializations of families of tight contact structures on $\R^3$, generalizing a classical result of Eliashberg \cite{Eliashberg}.

\textbf{Notation and conventions.}
Throughout this manuscript, $M$ will denote a compact, oriented $3$-manifold, possibly with non-empty boundary. Every contact structure $\xi$ on $M$ is assumed to be cooriented and positive. For a subset $C\subseteq M$, we denote by $\Op(C)$ an arbitrarily small and unspecified open neighborhood of $C$.

The group of orientation-preserving diffeomorphisms of $M$ that fix $\Op(\partial M)$ pointwise will be denoted by $\Diff(M)$, and $\Diff_0(M)$ will denote the path-connected component of the identity. When $M$ is equipped with a contact structure $\xi$, we denote by $\Cont(M,\xi)\subseteq \Diff(M)$ the subgroup of contactomorphisms, while $\Cont_0(M,\xi)=\Diff_0(M)\cap\Cont(M,\xi)$ denotes the subgroup of contactomorphisms smoothly isotopic to the identity. We denote by $\CStr(M)$ the space of contact structures on $M$ that coincide with an implicitly fixed contact structure on $\Op(\partial M)$, and by $\CStr(M,\xi)$ the path-connected component containing $\xi$.

The space of smooth embeddings $S\rightarrow M$ will always be denoted by $\Emb(S,M)$. Unless otherwise specified, our embeddings will always coincide on $\Op(\partial S)$ with a fixed reference embedding. The reference embedding will be specified when relevant. The only exception to this convention concerns embeddings of the $3$-ball: in $\Emb(\D^3,M)$, the boundary is allowed to move.

For an embedding $e\in \Emb(S,M)$, we denote by $\Emb(S,M,e)\subseteq \Emb(S,M)$ the path-connected component containing $e$. When $S$ is a surface and $\mathcal{F}$ is a $1$-dimensional singular foliation on $S$, we denote by $\Emb((S,\mathcal{F}),(M,\xi))\subseteq \Emb(S,M) $ the space of embeddings $S\rightarrow (M,\xi)$ whose induced characteristic foliation is precisely $\mathcal{F}$. The spaces of Legendrian and transverse embeddings $\NS^1\rightarrow (M,\xi)$ will be denoted by $\Leg(M,\xi)$ and $\Trans(M,\xi)$, respectively. For a Legendrian embedding $\Lambda$, we denote by $\Leg(M,\xi,\Lambda)\subseteq \Leg(M,\xi)$ the path-connected component containing $\Lambda$, and similarly for transverse embeddings.

For a subset $C$, the notation $\rel C$ will always mean relative to $\Op(C)$. Thus, $\Diff(M,\rel C)$ denotes the subgroup of diffeomorphisms that are the identity on $\Op(C)$, while $\Emb(S,M,\rel C)$ denotes the subspace of embeddings that coincide on $\Op(C)$ with the fixed reference embedding.

All the function spaces appearing in this article are equipped with the Whitney $C^\infty$ topology. One can show that all the relevant function spaces have the homotopy type of CW-complexes, so Whitehead's Theorem applies to upgrade all our weak homotopy equivalences to homotopy equivalences. Since in our applications we usually only need the weak homotopy type, we have chosen to work throughout with weak homotopy equivalences and leave this upgrade to the interested reader.

\textbf{Acknowledgements.}
We are grateful to R. Budney, V. Colin, Y. Eliashberg, F. Gironella, \'A. Gonz\'alez-Prieto, G. Mati\'c, H. Min, J. Muñoz-Ech\'aniz, D. Pancholi and \'A. del Pino for valuable conversations. The first author would like to thank R. Casals, Y. Eliashberg and V. Ginzburg for hosting him during a visit to California, where this project unexpectedly began to take shape. Part of this project was developed while he was visiting the Institute for Advanced Study in 2021, and he would like to acknowledge the hospitality of the Institute. The second author wants to thank Á. del Pino for his constant support. He would also like to thank the Geometry group at Utrecht University where he was given a nice environment to develop this and other projects. Finally, we thank Y. Eliashberg, D. Pancholi and G. S\'anchez for their encouragement.

The authors acknowledge support from the Spanish national grant with reference number PID2022-142024NB-I00(MINECO/FEDER) and the ICMAT–Severo Ochoa grant CEX2019-000904-S. The first author was supported by Beca de Personal Investigador en Formación UCM. During the development of this work, the second author was funded by Programa Predoctoral de Formación de Personal Investigador No Doctor
del Departamento de Educación del Gobierno Vasco.


\section{Basics in contact $3$-manifolds} \label{sec:Basics}


The section is organized as follows. First, we provide a brief review of convex surface theory; see \cite{HH19,EtnyreConvex,Honda,GirouxConvex,Massot} for further details. We then state the natural fibrations arising from Gray Stability and the contact isotopy extension theorems; see \cite{GeigesBook,GirouxMassot,MasFibrNotes}.

\subsection{Convex surface theory in contact $3$-manifolds}\label{sub:ConvexSurfacetheory} Let $(M,\xi)$ be a contact $3$-manifold. A properly embedded oriented surface $S\subseteq M$ is said to be \emph{convex} if there exists a contact vector field $X$ that is positively transverse to $S$. A \emph{convex embedding} of $S$ into $(M,\xi)$ is an embedding $e:S\rightarrow M$ such that $e(S)\subseteq M$ is convex. The $1$-dimensional singular foliation $\mathcal{F}:=e^*\xi=de^{-1}(\xi\cap de(TS))$ on $S$ is called the \emph{characteristic foliation} of the surface. A convex surface $S$ has a neighborhood in $(M,\xi)$ contactomorphic to $(S\times\R,\ker(f\,dt+\beta))$, where $t$ is the $\R$-coordinate, $\beta\in\Omega^1(S)$, and $f\in C^\infty(S)$. We call such a neighborhood \emph{vertically invariant}.

The \emph{dividing set} of $S$ is the embedded $1$-dimensional submanifold\[
\Gamma:=\{p\in S:X_p\in\xi_p\}.
\]
Although the dividing set depends on the contact vector field $X$, its isotopy class is well-defined because the space of contact vector fields positively transverse to $S$ is contractible. We will write $\Gamma_X$ instead of $\Gamma$ whenever we want to specify the contact vector field. The dividing set separates the surface into $S\setminus\Gamma_X=S_+\sqcup S_-$, where $S_+$ (resp. $S_-$) is the open subsurface on which $X$ is positively (resp. negatively) transverse to $\xi$.

A key feature of convex surfaces is that they are generic among smooth surfaces. 

\begin{theorem}[Giroux Approximation Theorem \cite{GirouxConvex,Honda,Kanda2}] \label{thm:GirApp} 
    Let $(M,\xi)$ be a contact $3$-manifold and let $S\subseteq M$ be a compact surface with Legendrian boundary. Then,
    \begin{itemize}
        \item If $\partial S=\emptyset$, then, after a $C^\infty$-small perturbation of $S$, we may assume that $S$ is convex.
        \item If $\partial S\neq\emptyset$ and the twisting of the contact planes along $\partial S$ with respect to $S$ is non-positive, then, after a $C^0$-small perturbation near $\partial S$ fixing $\partial S$, followed by a $C^\infty$-small perturbation elsewhere, we may assume that $S$ is convex.
    \end{itemize}
\end{theorem}

We will apply the previous theorem to surfaces that are already convex near their Legendrian boundary. Thus, in those applications, the perturbation can be chosen to be $C^\infty$-small.

Another fundamental result in convex surface theory is that tightness in a neighborhood of a convex surface can be characterized in terms of the topology of its dividing set.

\begin{theorem}[Giroux Tightness Criterion \cite{GirouxCriterion}]\label{thm:GirouxTightnessCriterion} Let $S\subseteq (M,\xi)$ be a convex surface in a contact $3$-manifold. Then, $S$ has a tight neighborhood if and only if 
\begin{itemize} 
\item $S$ is a sphere with connected dividing set, or 
\item $S$ is not a sphere and its dividing set does not contain any homotopically trivial curve. 
\end{itemize} 
\end{theorem}

This theorem illustrates a second fundamental feature of convex surfaces: the dividing set, rather than the precise characteristic foliation, carries the essential information about the contact structure in a neighborhood of the surface. This principle is formalized by the \emph{Giroux Realization Theorem}. To state it, we introduce some notation. Let $e:S\rightarrow (M,\xi)$ be a convex embedding with characteristic foliation $\mathcal{F}$ and dividing set $\Gamma$. We assume that $\partial S$ is Legendrian. We denote by $\Emb((S,\Gamma),(M,\xi))$ the space of convex embeddings $j:S\rightarrow M$ with dividing set $\Gamma$.

\begin{theorem}[Giroux Realization Theorem \cite{Giroux}]\label{thm:RealizationTheorem}
    The natural inclusion \[ \Emb((S,\mathcal{F}),(M,\xi))\hookrightarrow \Emb((S,\Gamma),(M,\xi))\]
    is a weak homotopy equivalence.
\end{theorem}
The deformations in the previous theorem are realized by \emph{admissible} isotopies, i.e. graphical deformations of $S\subseteq S\times \R$ inside a vertically invariant neighborhood. An important consequence of the Realization Theorem is the so-called Legendrian Realization Principle, due to Honda \cite{Honda}, which provides sufficient conditions for a collection of curves and arcs in a convex surface to become Legendrian by an admissible isotopy, making them part of the characteristic foliation.

An important result of Giroux is that the Realization Theorem may also be applied at the boundary of a contact $3$-manifold when studying the weak homotopy type of its contactomorphism group and, via the Gray fibration, the corresponding spaces of contact structures with prescribed boundary data:

\begin{proposition}\cite{Giroux,GirouxMassot}\label{prop:GirouxChangingBoundary}
Let $\CStr(M;\Gamma)$ be the space of contact structures on $M$ with convex boundary and dividing set $\Gamma\subseteq \partial M$. Then, the weak homotopy type of $\Cont(M,\xi)$, for $\xi\in\CStr(M;\Gamma)$, depends only on the connected component of $\CStr(M;\Gamma)$ containing $\xi$.
\end{proposition}

The following lemma states that the dividing sets of two convex surfaces alternate along a transverse Legendrian intersection.

\begin{lemma}\cite{Honda,Ka1997}\label{lem:DividingSetIntersection}
Let $S,S'\subseteq (M,\xi)$ be two embedded convex surfaces with dividing sets $\Gamma\subseteq S$ and $\Gamma'\subseteq S'$. Assume that $S$ and $S'$ intersect transversely along a Legendrian curve $\Lambda$. Then, between any two adjacent points of $\Gamma\cap\Lambda$ there is a point of $\Gamma'\cap\Lambda$, and vice versa.
\end{lemma}

A very useful feature of convex surfaces is that the formal invariants of a Legendrian knot can be easily computed from a convex Seifert surface:

\begin{lemma}[\cite{Kanda2}]\label{lem:FormalInvariantsDisk}
Let $\Lambda\subseteq (M,\xi)$ be an oriented Legendrian knot, and let $S$ be a convex Seifert surface for $\Lambda=\partial S$ with dividing set $\Gamma$. Then,
\begin{itemize}
    \item $\tb(\Lambda)=-\frac{1}{2}\#(\Gamma\cap\Lambda)$, and
    \item $\rot(\Lambda)=\chi(S_+)-\chi(S_-)$.
\end{itemize}
\end{lemma}

\subsection{Fibrations in contact topology}\label{sub:Fibrations} 

In the following lemma, we collect the key Serre fibrations that will be used throughout this manuscript.

\begin{lemma}[\cite{GirouxMassot,GeigesBook,MasFibrNotes}]\label{lem:fibrations} 
The following statements hold.
\begin{itemize} 
    \item [(i)] The map \[\Diff_0(M)\rightarrow\CStr(M,\xi), \qquad \varphi\mapsto \varphi_*\xi,\] is a fibration with fiber $\Cont_0(M,\xi)$.
    
    \item [(ii)] The map\[\Cont(M,\xi)\rightarrow\Leg(M,\xi),\qquad \varphi\mapsto\varphi\circ \Lambda,\] is a fibration with fiber $\Cont(C(\Lambda,M),\xi)$.
    
    \item [(iii)] The map \[ \Cont(M,\xi;\rel e(\partial S))\rightarrow\Emb((S,\mathcal{F}),(M,\xi)), \qquad \varphi\mapsto \varphi\circ e,\] is a fibration with fiber $\Cont(M,\xi;\rel e(S))$.
\end{itemize}
\end{lemma}

In (ii), we abuse notation by identifying the contactomorphism group of the Legendrian exterior $\Cont(C(\Lambda,M),\xi)$ with the subgroup $\Cont(M,\xi;\rel\Lambda)$ of contactomorphisms of $(M,\xi)$ that are the identity near $\Lambda$.


\section{Hatcher's Theorems}\label{sec:Hatcher}


In this section, we recall some results of Hatcher that will be used later in the article.

\subsection{The Smale Conjecture} We begin with Hatcher's affirmative answer to the \emph{Smale Conjecture}.

\begin{theorem}[Hatcher \cite{HatcherSmale}]\label{thm:SmaleConjecture}
The group of diffeomorphisms $\Diff(\D^3)$ of the $3$-ball fixing the boundary pointwise is weakly contractible.
\end{theorem} 
\begin{remark}
Recall that this result is equivalent to the usual formulation of the $3$-dimensional Smale conjecture, which asserts that the inclusion \[ \SO(4)\hookrightarrow \Diff(\NS^3) \] is a weak homotopy equivalence.
\end{remark}

\subsection{Embeddings of disks into irreducible $3$-manifolds} Let $M$ be an orientable, compact, connected, irreducible $3$-manifold with boundary. Let $i:\D^2\hookrightarrow M$ be a proper embedding and consider the space $\Emb(\D^2,M)$ of embeddings $\phi:\D^2\rightarrow M$ such that $\phi_{|\Op(\partial\D^2)}=i_{|\Op(\partial\D^2)}$.

\begin{theorem}[Hatcher \cite{HatcherHomeo,HatcherSurfaces}]\label{thm:DiskIrreducibleManifold}
The space $\Emb(\D^2,M)$ is weakly contractible.
\end{theorem}

This theorem, together with Theorem \ref{thm:SmaleConjecture}, allows us to determine the weak homotopy type of the diffeomorphism group of any handlebody. We include the argument since we will later replicate it in the contact setting for (multi)standard contact handlebodies.

\begin{corollary}\label{cor:DiffHandlebody}
Let $H^g$ be a genus $g$ handlebody. The diffeomorphism group $\Diff(H^g)$ is weakly contractible.
\end{corollary} 
\begin{proof}
Since $\Diff(\D^3)$ is weakly contractible by Theorem \ref{thm:SmaleConjecture}, it is enough to check that $\Diff(H^g)$ is weakly homotopy equivalent to $\Diff(H^{g-1})$ for every $g\geq1$. Let $e:\D^2\rightarrow H^g$ be a proper embedding of a non-separating disk. Postcomposition with $e$ defines a fibration \[ \Diff(H^g,\rel e(\D^2))\longrightarrow\Diff(H^g)\longrightarrow\Emb(\D^2,H^g). \] The fiber can be identified with $\Diff(H^{g-1})$, while the base is weakly contractible by Theorem \ref{thm:DiskIrreducibleManifold}. The result follows.
\end{proof}

The previous result also determines the weak homotopy type of the space of long unknots in the $3$-sphere. Our proof of Theorem \ref{thm:UnknotsS3Intro}(i) will later follow the same pattern in the contact setting.

\begin{theorem}[Hatcher \cite{HatcherSmale}]\label{thm:SmoothLongUnknots}
	The path-connected component $\Emb_{p,v}(\NS^1,\NS^3,U)$ of the space of smooth long embeddings containing the unknot is weakly contractible.
\end{theorem}
\begin{proof}
	Consider the fibration \[ \Diff(C(U,\NS^3)) \hookrightarrow \Diff(\NS^3,\rel p) \rightarrow \Emb_{p,v}(\NS^1,\NS^3,U). \]
    The total space $\Diff(\NS^3,\rel p)\cong \Diff(\D^3)$ is weakly contractible by Theorem \ref{thm:SmaleConjecture}, while the fiber $\Diff(C(U,\NS^3)) \cong \Diff(\D^2\times\NS^1)$ is weakly contractible by  Corollary \ref{cor:DiffHandlebody}.
\end{proof}
\begin{remark} In general, each path-connected component of the space of long embeddings in $\NS^3$ is a $K(\pi,1)$ space \cite{HatcherSmale}.
\end{remark}


\section{Eliashberg--Mishachev's Theorem}\label{sec:Eliash}


In this section, we recall the result of Eliashberg and Mishachev concerning tight contact structures on the $3$-disk \cite{EliashbergMishachevTight} and discuss some of its immediate consequences. As mentioned in the Introduction, our strategy for studying spaces of contact structures and contactomorphism groups is to cut a contact $3$-manifold along surfaces in order to reduce its complexity. Since the $3$-disk is the simplest possible piece arising from this procedure, the work of Eliashberg and Mishachev provides a foundational ingredient for the present article.

\subsection{Tight contact structures on the $3$-disk and $3$-sphere}

The starting point in the classification program for tight contact structures is the following result:

\begin{theorem}[Eliashberg \cite{Eliashberg}]\label{thm:EliashTightD3}
Any two tight contact structures on $\D^3$ that induce the same characteristic foliation along $\partial\D^3$ are isotopic relative to the boundary.
\end{theorem}

\begin{remark}
Observe that, by Darboux's Theorem, this implies the uniqueness, up to isotopy, of tight contact structures on $\NS^3$.
\end{remark}

Theorem \ref{thm:EliashTightD3} admits a parametric version. This was stated without proof in \cite{Eliashberg}; a full proof, due to Eliashberg and Mishachev, appeared recently in \cite{EliashbergMishachevTight}. There is also an approach to proving this result using techniques from convex surface theory in the thesis of D. Jänichen \cite{Janichen}.

\begin{theorem}[Eliashberg--Mishachev \cite{EliashbergMishachevTight}]\label{thm:ContactStructuresR3StandardInfinity}
The space $\CStr(\D^3,\xi_\std)$ of tight contact structures on $\D^3$ that coincide with the standard one near the boundary is contractible.
\end{theorem}

It follows from Giroux's Theorems \ref{thm:GirApp}, \ref{thm:RealizationTheorem} and \ref{thm:GirouxTightnessCriterion} that the same statement holds for any prescribed germ of a tight contact structure near the boundary:

\begin{corollary}\label{cor:ContactStructuresD3}
Let $(\D^3,\xi)$ be a tight ball. Then the space $\CStr(\D^3,\xi)$ of tight contact structures on $\D^3$ that coincide with $\xi$ near $\partial\D^3$ is contractible.
\end{corollary}

We now turn to the closed case of $\NS^3$. Consider the evaluation map
\begin{equation}\label{eq:EvaluationS3}
\ev_N:\CStr(\NS^3,\xi_\std)\longrightarrow \SO(4)/\U(2)=\NS^2
\end{equation}
defined by sending a tight contact structure on $\NS^3$ to its positively oriented unit normal vector at the north pole $N\in\NS^3$. This is a fibration whose fiber can be identified with $\CStr(\D^3,\xi_\std)$. Therefore, Theorem \ref{thm:ContactStructuresR3StandardInfinity} implies the following.

\begin{corollary}[Eliashberg--Mishachev]\label{cor:ContactStructuresS3}
The evaluation map (\ref{eq:EvaluationS3}) is a weak homotopy equivalence.
\end{corollary}

\begin{remark}
It is not difficult to describe a homotopy inverse of (\ref{eq:EvaluationS3}). Identify the space of isometric linear complex structures on $\C^2$ that induce the standard orientation with the $2$-sphere
\[
\NS^2=\NS^2(i,j,k)
=
\{J_{p_1,p_2,p_3}=p_1i+p_2j+p_3k:(p_1,p_2,p_3)\in\NS^2\}.
\]
Then a homotopy inverse of the evaluation map is given by the inclusion
\begin{equation}\label{eq:ComplexToTight}
i:\NS^2(i,j,k)\longrightarrow \CStr(\NS^3,\xi_\std),
\qquad
J\longmapsto \xi_J:=T\NS^3\cap JT\NS^3.
\end{equation}
\end{remark}

\subsection{Contactomorphisms of the tight $3$-disk and $3$-sphere}

The Gray fibration, together with Eliashberg--Mishachev's theorem and Hatcher's proof of the Smale conjecture (Theorem \ref{thm:SmaleConjecture}), determines the weak homotopy type of the contactomorphism groups of the tight $3$-disk and $3$-sphere.

We first state the result for the $3$-disk.

\begin{corollary}[Eliashberg--Mishachev]\label{coro:Contcontract}
The group $\Cont(\D^3,\xi_\std)$ of compactly supported contactomorphisms of the Darboux ball is contractible.
\end{corollary}
\begin{proof}
The result follows from the Gray fibration
\[
\Cont(\D^3,\xi_\std)\hookrightarrow
\Diff(\D^3)\longrightarrow
\CStr(\D^3,\xi_\std),
\]
together with Theorems \ref{thm:SmaleConjecture} and \ref{thm:ContactStructuresR3StandardInfinity}.
\end{proof}

As for contact structures, the same statement holds for any other tight boundary condition.

\begin{corollary}\label{coro:ContactomorphismsTightBall}
Let $(\D^3,\xi)$ be a tight $3$-disk. The group $\Cont(\D^3,\xi)$ of compactly supported contactomorphisms of $(\D^3,\xi)$ is contractible.
\end{corollary}


Finally, we turn to the case of the $3$-sphere.

\begin{corollary}\label{cor:ContsctomorphismSphere}
The natural inclusion
\begin{equation}\label{eq:InclusionContU(2)}
i:\U(2)\hookrightarrow\Cont(\NS^3,\xi_\std)
\end{equation}
is a weak homotopy equivalence.
\end{corollary}
\begin{proof}
Consider the following commuting diagram, whose rows are fibrations:
\begin{displaymath}
\xymatrix@M=10pt{
\U(2) \ar@{^{(}->}[d] \ar[r]&
\SO(4) \ar@{^{(}->}[d] \ar[r]&
\NS^2(i,j,k)=\SO(4)/\U(2) \ar@{^{(}->}[d] \\
\Cont(\NS^3,\xi_\std) \ar[r]&
\Diff(\NS^3) \ar[r]&
\CStr(\NS^3,\xi_\std).
}
\end{displaymath}
The result follows from Theorem \ref{thm:SmaleConjecture} and Corollary \ref{cor:ContactStructuresS3}.
\end{proof}


\section{The Embedding Microfibration Method}\label{sec:MicrofibrationTrick}

In this section we introduce the \emph{Embedding Microfibration Method}, which can be viewed as a parametric variation of Colin's \emph{Discretization Technique} \cite{Colin}. We will use this idea to prove several results concerning spaces of convex embeddings of disks with Legendrian boundary and spheres into a tight contact $3$-manifold. In general, this method can be applied when the following two conditions are satisfied:
\begin{itemize}
    \item [(i)] The space of convex embeddings $S\hookrightarrow(M,\xi)$ with fixed characteristic foliation $\mathcal{F}$ is dense in the space of smooth embeddings.
    \item [(ii)] For a vertically invariant neighborhood $(S\times \R,\xi)$, we understand, at the homotopy level, the inclusion of the space of convex embeddings $S\hookrightarrow (S\times \R,\xi)$ with fixed characteristic foliation $\mathcal{F}$ into the corresponding space of smooth embeddings.
\end{itemize}
Colin's technique \cite{Colin} requires precisely the same hypotheses at the level of $\pi_0$.

\subsection{The Microfibration Lemma}\label{sub:Microfibration}

The method is based on the notion of a \emph{Serre microfibration}, introduced by M. Gromov in \cite{GromovPartial}. The difference between Serre fibrations and Serre microfibrations is that, in the latter case, the lifting property is only guaranteed for sufficiently small times. The precise definition is as follows.

\begin{definition}
A map $p:E\rightarrow B$ is a \em Serre microfibration \em if for every $k\geq 0$ and pair of maps $h:\D^k\times[0,1]\rightarrow B$ and $g:\D^k\times\{0\}\rightarrow E$ such that $p\circ g=h_{|\D^k\times\{0\}}$ there exists a positive number $\varepsilon>0$ and a map $\hat{g}:\D^k\times[0,\varepsilon]\rightarrow E$ such that
\begin{itemize}
    \item $\hat{g}_{|\D^k\times\{0\}}=g$ and
    \item $p\circ \hat{g}=h_{|\D^k\times[0,\varepsilon]}$.
\end{itemize}
\end{definition}

The following result, also hinted by M. Gromov \cite{GromovPartial}, gives a sufficient condition for a microfibration to be a fibration.

\begin{lemma}[Microfibration Lemma \cite{Weiss}] \label{lem:micro}
Let $p:E\rightarrow B$ be a Serre microfibration with weakly contractible non-empty fibers. Then $p$ is a Serre fibration.
\end{lemma}

\subsection{The method}
Let $\Emb(N,M)$ be the space of smooth embeddings from a compact manifold $N$ into a target manifold $M$. As usual, if $N$ has non-empty boundary we will always require every embedding to be fixed on an open neighborhood of the boundary. For a given embedding $\phi\in\Emb(N,M)$  we consider the family of open neighborhoods
$$ \mathcal{U}_{\varphi}(\varepsilon)= \{ \psi\in \Emb(N,M): \forall x\in N, d(\psi(x), \phi(N))< \varepsilon \},$$
for each $\varepsilon>0$.

Let $E_0\subseteq E_1\subseteq \Emb(N,M)$ be two subspaces of $\Emb(N,M)$.

\begin{definition}
The inclusion $j:E_0 \hookrightarrow E_1$ is a \em $C^0$-dense homotopy equivalence \em if: 
\begin{itemize}
    \item [(i)] $j$ is a homotopy equivalence.
    \item [(ii)] For $\varepsilon>0$ and every family of embeddings $\varphi^{k}\in E_1$, $k\in K$, where $K$ is a compact CW-complex, such that $\varphi^k\in E_0$ for $k\in G\subseteq K$, where $G$ is a subcomplex of $K$, there exists a homotopy $\varphi^{k}_{t}\in E_1$, $t\in[0,1]$, such that 
    \begin{itemize}
    \item [(a)] $\varphi^{k}_0=\varphi^{k}$ for $k\in K$,
    \item [(b)] $\varphi^{k}_t=\varphi^{k}$ for $k\in G$,
    \item [(c)] $\varphi^{k}_1\in E_0$, for $k\in K$, and 
    \item [(d)] For every $k\in K$, the homotopy $\varphi^{k}_t\in \mathcal{U}_{\varphi^k}(\varepsilon)$, for every $t\in[0,1]$.
    \end{itemize}
\end{itemize}
\end{definition}

If we drop condition (i) from the definition we will say that the inclusion is a \em $C^0$-dense weak homotopy equivalence. \em In our examples, the $C^0$-dense weak homotopy equivalence property will imply the $C^0$-dense homotopy equivalence property because of Whitehead's Theorem. 

We can now state our embedding microfibration method, which provides a criterion for verifying the $C^0$-dense weak homotopy equivalence property. 

\begin{theorem} \label{thm:naked}
Let $E_0\subseteq E_1\subseteq \Emb(N,M)$ be two subspaces such that:
\begin{itemize}
\item [(i)] $E_0$ is dense in $E_1$.
\item [(ii)] For every embedding $\varphi\in E_1$, there exists $\varepsilon(\varphi)>0$ such that 
the natural map \[E_0 \cap \mathcal{U}_{\varphi}(\varepsilon)\to E_1\cap\mathcal{U}_{\varphi}(\varepsilon)\] induces a weak homotopy equivalence for every $0<\varepsilon<\varepsilon(\phi)$.
\end{itemize}
Then, the inclusion $j:E_0\hookrightarrow E_1$ is a $C^0$-dense weak homotopy equivalence.
\end{theorem}
\begin{remark}
Observe that condition (i) is formally implied by (ii). We state them separately since, in our applications, both conditions arise from different geometric arguments. 
\end{remark}
\begin{proof}
Let $(K,G)$ be a compact CW-pair and $\varphi^k\in E_1$, $k\in K$, be a family of embeddings such that $\varphi^k\in E_0$ for every $k\in G$. First, choose $\varepsilon>0$ sufficiently small so that all the maps 
\[ E_0 \cap \mathcal{U}_{\varphi^k}(\varepsilon)\to E_1\cap \mathcal{U}_{\varphi^k}(\varepsilon) \]
are weak homotopy equivalences. This can be done by a standard compactness argument.  

We denote the homotopy fiber of the inclusion $j:E_0\hookrightarrow E_1$ over $\varphi\in E_1$ by 
\[ F_\varphi=\{\varphi_t\in E_1:t\in[0,1],\varphi_0=\varphi, \varphi_1\in E_0\}. \] The key observation is that the subspace of $F_{\varphi^k}$ consisting of paths that remain sufficiently close to $\varphi^k$ is non-empty and weakly contractible by (i) and (ii). More precisely, the space \[F_{\varphi^k}(\varepsilon):=\{\varphi^k_t\in F_{\varphi^k}: \varphi^k_t\in U_{\varphi^k}(\varepsilon) \text{ for every } t\in [0,1]\}\] 
is non-empty and weakly contractible. These spaces will be the fibers of a suitable microfibration, which we now define. 

Let \[ X_\varepsilon:= \{ (k,\varphi^k_t): k\in K, \varphi^k_t\in F_{\varphi^k}(\varepsilon)\}.  \] Projection onto the first component defines a microfibration \[ F_{\varphi^k}(\varepsilon) \hookrightarrow X_\varepsilon \rightarrow K. \]  The microfibration property follows from the triangle inequality: for $\delta>0$ small enough, if $d_{C^0}(\varphi_0, \varphi_1)< \delta$ then we have that  $\mathcal{U}_{\varphi_1}(\varepsilon-\delta) \subset \mathcal{U}_{\varphi_0}(\varepsilon)$. 

Since the fiber is non-empty and weakly contractible, Lemma \ref{lem:micro} implies that the microfibration is, in fact, a Serre fibration. 

To conclude, observe that the family $\varphi^{k}_t= \varphi^k$, $(k,t)\in G\times [0,1]$, defines a section of this fibration over $G\subseteq K$. Since the fiber is weakly contractible this section extends to a section $\varphi^k_t$ over $K$, solving the problem.
\end{proof}


\section{Spaces of convex disks and spheres in a tight contact $3$-manifold}\label{sec:ConvexSurfaces}


The purpose of this section is to apply the microfibration method (Theorem \ref{thm:naked}) to several spaces of convex disks and spheres in tight contact $3$-manifolds. Tightness will be crucial in establishing the density property required to apply the method.

\subsection{Standard disks in a tight contact $3$-manifold}\label{subsubsec:StandardDisks}

\begin{figure}[!h]
  \centering
{\includegraphics[width=0.35\textwidth]{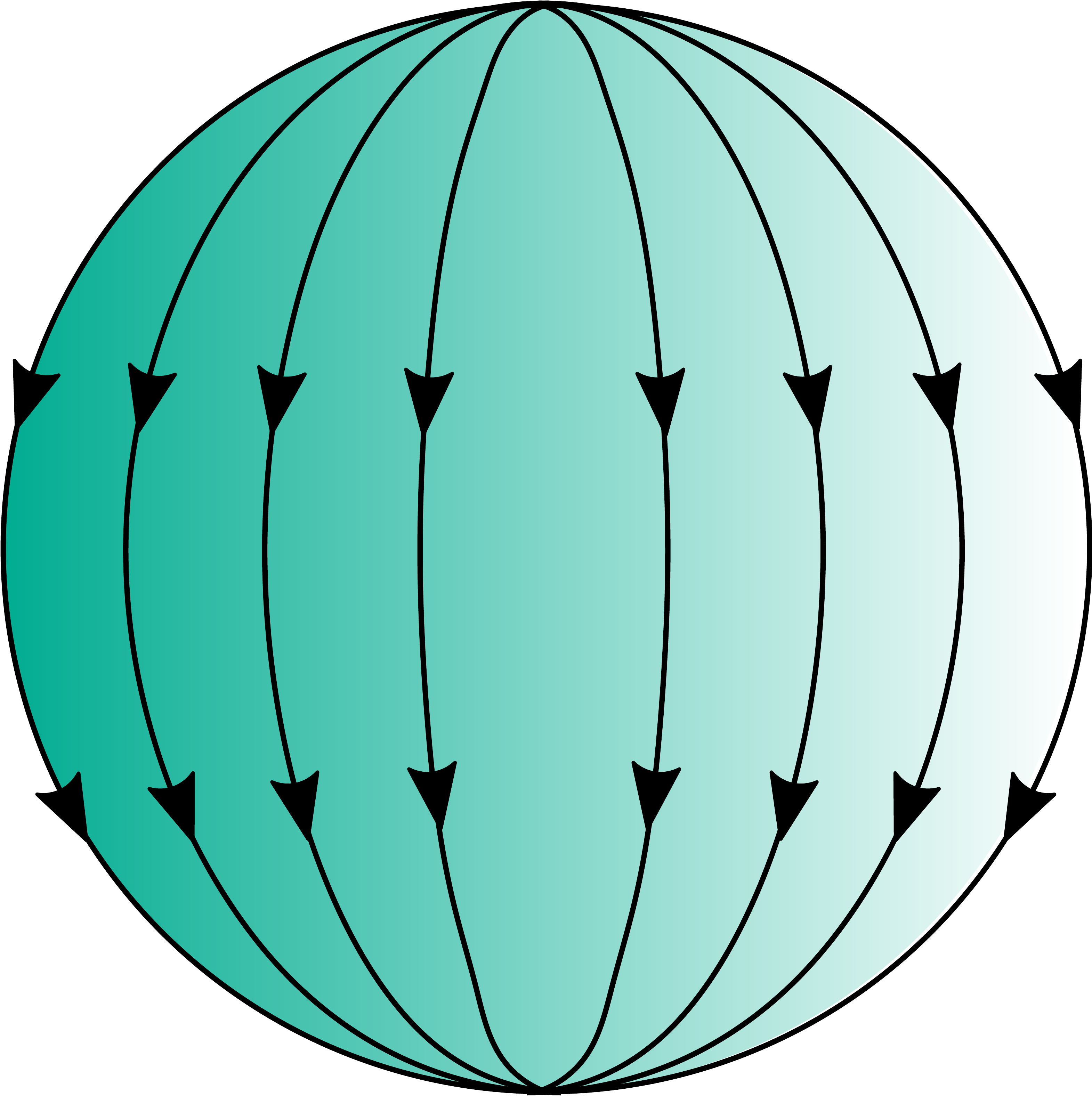}}
  \caption{The standard convex disk.\label{fig:standard}}
\end{figure}

Denote by $(\D^2, \mathcal{F}_\std)$ the standard foliation of the disk with boundary the Legendrian unknot with $\tb=-1$. That is, the characteristic foliation has two elliptic points of opposite signs at the boundary and Legendrian leaves joining them. See Figure \ref{fig:standard}. We will say that an embedding $e:\D^2\rightarrow(M,\xi)$ of a disk into a contact $3$-manifold is \emph{standard} if it is convex and $(\D^2,e^*\xi)=(\D^2,\mathcal{F}_\std)$. In particular, the boundary of $e(\D^2)$ is a standard Legendrian unknot, that is, with $\tb=-1$.

Let $e:\D^2\rightarrow(M,\xi)$ be a standard embedding of a disk into a contact $3$-manifold. Denote by $\Emb(\D^2,M)$ the space of smooth embeddings $f:\D^2\rightarrow M$ such that $f_{|\Op(\partial\D^2)}=e_{|\Op(\partial\D^2)}$. Similarly, define $\Emb((\D^2, \mathcal{F}_\std),(M,\xi))$ as the subspace of $\Emb(\D^2,M)$ formed by embeddings that are standard.

Our main result reads as follows.
\begin{theorem} \label{thm:standardDisks}
  Let $(M,\xi)$ be a tight contact $3$-manifold. The inclusion \[\Emb((\D^2, \mathcal{F}_\std),(M,\xi))\hookrightarrow\Emb(\D^2,M)\]is a $C^0$-dense weak homotopy equivalence.
\end{theorem}
\begin{proof}
We will apply Theorem \ref{thm:naked} with $E_0=\Emb((\D^2, \mathcal{F}_\std),(M,\xi))$ and $E_1=\Emb(\D^2,M)$. We need to check two things:
\begin{enumerate}
\item [(i)] The space of standard disks is $C^0$-dense inside the space of smooth disks. Since $(M,\xi)$ is tight and the boundary of the disk is a Legendrian unknot with $\tb=-1$, there is a unique configuration of dividing set for convex disks with Legendrian boundary with $\tb=-1$, namely a diameter. This follows from Giroux's Tightness Criterion (Theorem \ref{thm:GirouxTightnessCriterion}) and Lemma \ref{lem:FormalInvariantsDisk}. Therefore, the density property follows from Giroux's approximation and realization theorems, i.e. Theorems \ref{thm:GirApp} and \ref{thm:RealizationTheorem}.
\item [(ii)] For every embedding $\varphi\in \Emb(\D^2,M)$, the inclusion  $E_0 \cap U_{\varphi}(\varepsilon)\hookrightarrow U_{\varphi}(\varepsilon)\cap E_1=U_{\varphi}(\varepsilon)$ is a weak homotopy equivalence for every $\varepsilon>0$ small enough.  
\end{enumerate}
Observe that if $\varepsilon>0$ is small then $U_{\varphi}(\varepsilon)$ coincides with the space of smooth disk embeddings in a tight ball. Therefore, property (ii), and hence the result, follows from Lemma \ref{lem:DisksBall} proved below.
\end{proof}

\begin{lemma}\label{lem:DisksBall}
	Let $(\D^3,\xi)$ be a tight contact $3$--ball. Fix coordinates $(x,y,z)\in\D^3$ and assume that $e:\D^2(x,y)\hookrightarrow(\D^3,\xi), (x,y)\mapsto (x,y,0)$ is a standard embedding. Then, 
	the inclusion \[\Emb((\D^2, \mathcal{F}_\std),(\D^3,\xi))\hookrightarrow\Emb(\D^2, \D^3)\] is a weak homotopy equivalence. In particular, the space $\Emb((\D^2, \mathcal{F}_\std),(\D^3,\xi))$ is contractible.
\end{lemma}
\begin{proof}
Consider the $3$-balls $\D^3_1=\{(x,y,z)\in\D^3:z\geq0\}$ and $\D^3_2=\{(x,y,z)\in\D^3:z\leq0\}$. There is a natural commuting diagram in which the rows are fibrations \begin{displaymath}
\xymatrix@M=10pt{
			\Cont(\D^3_1, \xi)\times\Cont(\D^3_2,\xi)\ar@{^{(}->}[r]\ar@{^{(}->}[d]& \Cont(\D^3,\xi)\ar[r]\ar@{^{(}->}[d]& \Emb((\D^2, \mathcal{F}_\std),(\D^3,\xi)) \ar@{^{(}->}[d]\\
		\Diff(\D^3_1)\times\Diff(\D^3_2) \ar@{^{(}->}[r]& \Diff(\D^3) \ar[r]& \Emb(\D^2,\D^3). }
	\end{displaymath}
The statement follows from the Five Lemma, Theorem \ref{thm:SmaleConjecture} and Corollary \ref{coro:ContactomorphismsTightBall}.
\end{proof}

\subsection{Multistandard disks in a tight contact $3$-manifold}

Consider a convex embedding $e:\D^2\rightarrow (M,\xi)$ into a tight contact $3$-manifold, with characteristic foliation $\mathcal{F}$ and Legendrian boundary $\Lambda=e(\partial \D^2)$. The argument in the proof of Theorem \ref{thm:standardDisks} applies to the space \[E_0=\Emb((\D^2,\mathcal{F}),(M,\xi))\] of convex embeddings of disks with fixed Legendrian boundary $\Lambda$ and fixed characteristic foliation $\mathcal{F}$, provided that $E_0$ is dense in the corresponding space of smooth disk embeddings. Indeed, Lemma \ref{lem:DisksBall} requires no assumption on the characteristic foliation.

The density property will be guaranteed if, once we fix the Legendrian boundary $\Lambda=e(\partial \D^2)$, the dividing set is unique up to isotopy. 

\begin{figure}[!h]
  \centering
{\includegraphics[width=0.35\textwidth]{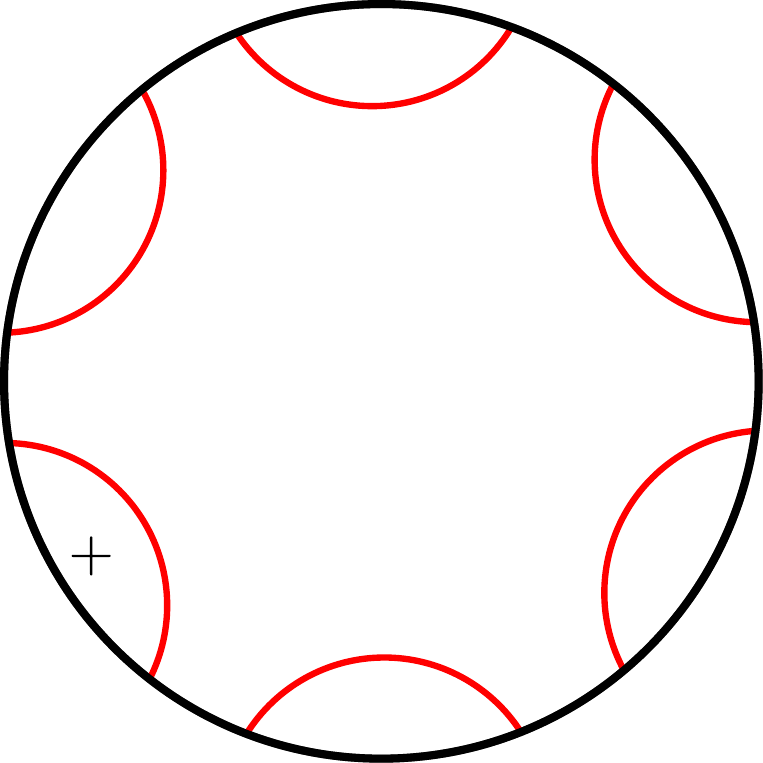}}
  \caption{The unique configuration of the dividing set for a disk with $\tb(\Lambda)=-6$ and $\rot(\Lambda)=5$.\label{fig:MaximalRotation}}
\end{figure}

A sufficient numerical condition ensuring the uniqueness of the dividing set when $(M,\xi)$ is tight is \begin{equation}\label{eq:MaxRot} \tb(\Lambda)+|\rot(\Lambda)|=-1. \end{equation} Indeed, as a consequence of Giroux's Criterion (Theorem \ref{thm:GirouxTightnessCriterion}) and Lemma \ref{lem:FormalInvariantsDisk}, the dividing set has $|\tb(\Lambda)|$ components, all of them $\partial$-parallel. Moreover, the sign of the region cobounded by each component of $\Gamma$ that does not contain any other component of $\Gamma$ agrees with the sign of $\rot(\Lambda)$. See Figure \ref{fig:MaximalRotation}. We will say that a convex disk $(\D^2,\mathcal{F})$ is \emph{multistandard} if its boundary Legendrian unknot satisfies (\ref{eq:MaxRot}), or, equivalently, its dividing set is $\partial$-parallel. 

We immediately obtain:

\begin{theorem}\label{thm:MultiStandardDisks}
    Let $(M,\xi)$ be a tight contact $3$-manifold and $(\D^2,\mathcal{F})$ a multistandard convex disk. Then, the inclusion $$\Emb((\D^2,\mathcal{F}),(M,\xi))\hookrightarrow\Emb(\D^2,M)$$ is a $C^0$-dense weak homotopy equivalence.
\end{theorem}

\begin{figure}[!h]
  \centering
{\includegraphics[width=1\textwidth]{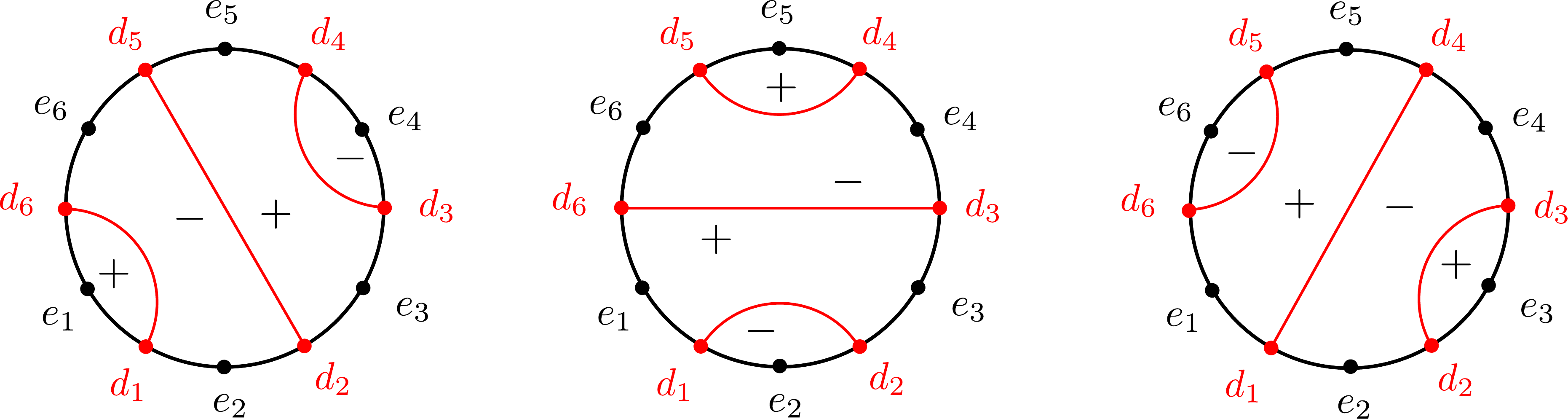}}
  \caption{The three possible configurations of the dividing set for a disk with $\tb(\Lambda)=-3$ and $\rot(\Lambda)=0$.\label{fig:TB-3}}
\end{figure}

\begin{proof}[Proof of Theorem \ref{thm:DisksIntro}]
The result follows by combining Theorems \ref{thm:standardDisks} and \ref{thm:MultiStandardDisks}, and Giroux's Realization Theorem \ref{thm:RealizationTheorem}.
\end{proof}

\subsection{Mini-disks in a tight contact $3$-manifold}
Let $(M,\xi)$ be a tight contact $3$-manifold and let \[E_0=\Emb((\D^2,\mathcal{F}_{mini}),(M,\xi))\] denote the space of convex disk embeddings with fixed transverse boundary $T$ and fixed characteristic foliation $\mathcal{F}_{mini}$. The characteristic foliation $\mathcal{F}_{mini}$ consists of a single positive elliptic singularity in the interior, with Legendrian leaves emanating from it. Such disks are referred to as \emph{mini-disks}, since they can be contracted to the elliptic point by a contact isotopy (without fixing the boundary). It follows that the boundary $T$ is a transverse unknot with self-linking number $\operatorname{sl}(T)=-1$. 

A consequence of Giroux's Elimination Lemma \cite{GirouxConvex} is that, when the ambient contact manifold is tight, $E_0$ is dense in the corresponding space of smooth disks; this  fact was used in \cite{EliashbergTransverse}. Thus, the first hypothesis in Theorem \ref{thm:naked} is satisfied. The second follows as in Lemma \ref{lem:DisksBall}. Therefore, we conclude:

\begin{theorem}\label{thm:MiniDisks}
Let $(M,\xi)$ be a tight contact $3$-manifold. Then, the inclusion 
\[ \Emb((\D^2,\mathcal{F}_{mini}),(M,\xi))\hookrightarrow \Emb(\D^2,M) \] is a $C^0$-dense weak homotopy equivalence.
\end{theorem}

\subsection{Standard spheres in a tight contact $3$-manifold}
In this final subsection, we provide the details of how to apply the microfibration method to the standard convex sphere. Let $\Emb((\NS^2,\mathcal{F}_\std),(M,\xi))$ denote the space of convex embeddings $e:\NS^2\rightarrow(M,\xi)$ that are \emph{standard}, meaning that the characteristic foliation $e^*\xi=\mathcal{F}_\std$ coincides with that induced on the boundary of a fixed Darboux ball. We also denote by $\Emb(\NS^2,M,e)$ the space of embeddings smoothly isotopic to $e$ and by
\[
\Emb((\NS^2,\mathcal{F}_\std),(M,\xi),e)\subseteq \Emb(\NS^2,M,e)
\]
the subspace consisting of standard embeddings. The following result can be regarded as a multiparametric version of Colin's result in \cite{Colin}.

\begin{theorem}\label{thm:Standard_spheres}
Let $(M,\xi)$ be a tight contact $3$-manifold and $e:\NS^2\rightarrow (M,\xi)$ be a standard embedding. Then, $$ \pi_k(\Emb(\NS^2,M,e),\Emb((\NS^2,\mathcal{F}_\std),(M,\xi),e))\cong\pi_k(\SO(3),\U(1)),$$ for any $k\geq1$. 
\end{theorem}

We will prove this result using the microfibration method of Theorem \ref{thm:naked}. We first establish some observations about spherical annuli that are needed to verify property (ii) of the method.

The differential topology statement about spherical annuli that we need is the following:

\begin{lemma} \label{lem:localEmb} The following statements hold.
	\begin{itemize}
		\item The diffeomorphism group $\Diff(\NS^2\times[-1,1])$ is weakly homotopy equivalent to $\Omega\SO(3)$. 
		\item Fix the inclusion $e:\NS^2\hookrightarrow\NS^2\times\{0\}\subseteq \NS^2\times[-1,1]$. The evaluation at the north pole of the $1$-jet defines a weak homotopy equivalence
		$$ D\ev_N:\Emb(\NS^2,\NS^2\times[-1,1],e)\rightarrow \SO(3), f\mapsto d_N f.$$
	\end{itemize}
\end{lemma}
\begin{proof}
	The first statement follows from the Smale Conjecture (Theorem \ref{thm:SmaleConjecture}) since $\Diff(\NS^2\times[-1,1])$ is the fiber of the fibration 
	$$ \Diff(\NS^2\times[-1,1])\hookrightarrow \Diff(\D^3)\rightarrow \Emb(\D^3(\varepsilon),\D^3).$$
	Indeed, the total space is contractible and the base is weakly homotopy equivalent to $\SO(3)$ by the Alexander trick. 
	
	Similarly, consider the commutative diagram 
		\begin{displaymath} 
	\xymatrix@M=10pt{
			\Omega\SO(3)\times\Omega\SO(3)\ar@{^{(}->}[r]\ar@{^{(}->}[d]& \Omega\SO(3)\ar[r]\ar@{^{(}->}[d]& \SO(3)\ar@{^{(}->}[d]\\
		\Diff(\NS^2\times[-1,0])\times\Diff(\NS^2\times[0,1]) \ar@{^{(}->}[r]& \Diff(\NS^2\times[-1,1]) \ar[r]& \Emb (\NS^2,\NS^2\times[-1,1],e), }
	\end{displaymath}
	where the rows are fibrations. It follows from the first statement and the Five Lemma that the inclusion $\SO(3)\hookrightarrow\Emb(\NS^2,\NS^2\times[-1,1],e)$ is a weak homotopy equivalence and, thus, the second statement.
\end{proof}

The contact analogue of the previous lemma also follows:

\begin{lemma} \label{lem:localSEmb}
	Let $(\NS^2\times[-1,1],\xi)$ be a tight spherical annulus. Assume that the inclusion $e:\NS^2\hookrightarrow\NS^2\times\{0\}\subseteq \NS^2\times[-1,1]$ is a standard embedding.
		\begin{itemize}
		\item The contactomorphism group $\Cont(\NS^2\times[-1,1],\xi)$ is weakly homotopy equivalent to $\Omega\U(1)$. 
		\item  The evaluation at the north pole of the $1$-jet defines a weak homotopy equivalence
		$$ D\ev_N:\Emb((\NS^2,\mathcal{F}_\std),(\NS^2\times[-1,1],\xi),e)\rightarrow \U(1), f\mapsto d_N f.$$
	\end{itemize}
\end{lemma}
\begin{proof}
	As usual, it follows from the results due to Giroux presented in Section \ref{sec:Basics} that it is enough to check the statement for the standard tight contact structure $\xi_\std$ on $\NS^2\times[-1,1]$, i.e. the canonical $\R$-invariant neighborhood of a standard convex sphere. 
	
	In the case of $(\NS^2\times[-1,1],\xi_\std)$ the statement follows in the same way as in the smooth case. Recall that the group $\Cont(\D^3,\xi_\std)$ is contractible and that the space $\CEmb((\D^3(\varepsilon),\xi_\std),(\D^3,\xi_\std))$ of Darboux ball embeddings is weakly homotopy equivalent to $\U(1)$ by the contact Alexander trick. 
\end{proof}
\begin{remark}
In the weak homotopy equivalence above, $\Omega\U(1)$ accounts for the powers of a contactomorphism first observed by Gompf \cite{Gompf}. This contactomorphism was termed the \emph{contact Dehn twist} along $\NS^2\times\{0\}$ in \cite{FME22}.
\end{remark}

\begin{proof}[Proof of Theorem \ref{thm:Standard_spheres}]
The $1$-jet evaluation of a smooth embedding $\varphi\in\Emb(\NS^2,M,e)$ at the north pole $N\in\NS^2$ defines a fibration $\Emb(\NS^2,M,e)\rightarrow \Fr^+(M)$. In the same way, for the standard case we have a fibration $\Emb((\NS^2,\mathcal{F}_\std),(M,\xi),e)\rightarrow \CFr(M,\xi)$. Both fibrations fit in a commuting diagram
	\begin{displaymath} 
	\xymatrix@M=10pt{
			F_0\ar@{^{(}->}[r]\ar@{^{(}->}[d]& 
\Emb((\NS^2,\mathcal{F}_\std),(M,\xi),e) \ar[r]\ar@{^{(}->}[d]& \CFr(M,\xi) \ar@{^{(}->}[d]\\
		F_1 \ar@{^{(}->}[r]& \Emb(\NS^2,M,e) \ar[r]& \Fr^+(M), }
	\end{displaymath}
 where the fiber $F_1$ is weakly homotopy equivalent to $E_1=\Emb(\NS^2,M;\rel N,e)\subseteq \Emb(\NS^2,M,e)$, the subspace of smooth embeddings which coincide with $e$ over $\Op(N)$, and $F_0$ is weakly homotopy equivalent to $E_0=\Emb((\NS^2,\mathcal{F}_\std),(M,\xi),\rel N,e)=E_1\cap \Emb((\NS^2,\mathcal{F}_\std),(M,\xi),e)$.
 
 Therefore, it is enough to check that the inclusion $$j:E_0\hookrightarrow E_1$$ is a weak homotopy equivalence. We will use Theorem \ref{thm:naked} to do it. As usual, the density property (i) follows easily from the tightness of $(M,\xi)$ and Giroux's convex surface theory. Observe that if $\varphi\in E_0$ is an embedding and $\varepsilon>0$ is small enough then $\mathcal{U}_\varphi (\varepsilon)$ can be identified with the space of smooth embeddings $\NS^2\hookrightarrow \NS^2\times[-1,1]$ which coincide with $e$ over $\Op(N)$. It follows that property (ii) in the hypothesis of Theorem \ref{thm:naked} can be deduced easily from Lemmas \ref{lem:localEmb} and \ref{lem:localSEmb}. 
 \end{proof}

\begin{proof}[Proof of Theorem \ref{thm:SpheresIntro}]
    By Giroux's Tightness Criterion and Realization Theorems \ref{thm:GirouxTightnessCriterion} and \ref{thm:RealizationTheorem}, the result follows immediately from Theorem \ref{thm:Standard_spheres}.
\end{proof}


\section{Contact handle attachments}\label{sec:Gluing}


In this section, we define our contact handle attachment maps and state and prove the technical version of Theorem \ref{thm:ContactHandleIntro}, our main gluing theorem for families of tight contact structures. We then turn to applications of this result, including Theorem \ref{thm:StandardHandlebodyINTRO}, Theorem \ref{thm:UnknotsS3Intro}(i), and several observations about contact cell decompositions.

\subsection{Contact handle attachment maps}

Let $(M,\xi)$ be a tight contact $3$-manifold with convex boundary. We will now define the contact handle attachment maps with domain $\CStr(M,\xi)$.

The effect of a contact $0$-handle attachment is to take the disjoint union with a Darboux ball $(H^0=\D^3,\xi_\std)$, producing the tight contact manifold $(M\sqcup H^0,\xi_{H^0})$. The contact $0$-handle attachment map is simply \begin{equation}\label{eq:Contact0Handle}
\mathcal{CH}_0:\CStr(M,\xi)\longrightarrow \CStr(M\sqcup H^0,\xi\sqcup\xi_\std), \qquad \hat{\xi}\mapsto \hat{\xi}_{H^0}.
\end{equation}

The contact $1$-handle attachment map will be defined under the assumption that $(M,\xi)$ contains two disjoint embedded standard disks $D_0,D_1\subseteq \partial M$ whose orientations induced from $\partial M$ are opposite. We can then build a new contact manifold by identifying $D_0$ and $D_1$ via an $I$-invariant neighborhood of a standard disk $(H^1=\D^1\times \D^2,\xi_{\mathrm{inv}})$, and rounding corners as in \cite{Honda}. The resulting contact manifold will be denoted by $(M\cup H^1,\xi_{H^1})$. The contact $1$-handle attachment map is then defined as \begin{equation}\label{eq:Contact1Handle}
\mathcal{CH}_1:\CStr(M,\xi)\longrightarrow \CStr(M\cup H^1,\xi_{H^1}), \qquad \hat{\xi}\mapsto \hat{\xi}_{H^1}.
\end{equation}

\begin{remark}
    This is not the standard description of a contact $1$-handle attachment, but it is equivalent to it. See Appendix A of \cite{HH19} for a nice description of contact $1$-handles from which the equivalence follows easily.
\end{remark} 
We conclude with contact $3$-handles. For this, we assume that $(M,\xi)$ has a spherical component $S\subseteq \partial M$ which is a standard sphere. The effect of attaching a contact $3$-handle along $S$ is simply to cap off this sphere with a Darboux ball $(H^3=\D^3,\xi_\std)$, thereby producing the manifold $(M\cup H^3,\xi_{H^3})$. At the level of spaces, the contact $3$-handle attachment map is \begin{equation}\label{eq:Contact3Handle}
\mathcal{CH}_3:\CStr(M,\xi)\longrightarrow \CStr(M\cup H^3,\xi_{H^3}), \qquad \hat{\xi}\mapsto \hat{\xi}_{H^3}.
\end{equation} Observe that the image of $\mathcal{CH}_3$ is naturally contained in $\CStr(M\cup H^3,\xi_{H^3},\rel p)$, where $p\in S$, the subspace of tight contact structures fixed near a point of the sphere.

\subsection{The gluing theorem}

We now state and prove our main gluing result. We retain the notation from the previous subsection.

\begin{theorem}\label{thm:Gluing}
    The maps $\mathcal{CH}_0$ and $\mathcal{CH}_1$ are weak homotopy equivalences.

    The map $\mathcal{CH}_3$ is a weak homotopy equivalence onto $\CStr(M\cup H^3,\xi_{H^3},\rel p)$.
\end{theorem}
\begin{proof}
    The map $\mathcal{CH}_0$ is a weak homotopy equivalence by Eliashberg--Mishachev's Theorem \ref{thm:ContactStructuresR3StandardInfinity}. 

    We turn to the case of $1$-handles. The crucial observation is that if $D\subseteq H^1$ is the cocore disk, then a contact structure
$\hat{\xi}\in \CStr(M\cup H^1,\xi_{H^1})$ lies, up to a \emph{canonical} isotopy obtained by flowing along the normal directions to $D$, in the image of
$\mathcal{CH}_1$ precisely when $D$ is standard. We need to prove that this can be achieved in families.

    That is, if we denote by $\mathcal{G}\subseteq \CStr(M\cup H^1,\xi_{H^1})$ the subspace of contact structures for which the cocore disk $D$ is standard, then we must prove that the inclusion \[\mathcal{G}\hookrightarrow \CStr(M\cup H^1,\xi_{H^1}) \] is a weak homotopy equivalence.
    
     Let $(K,G)$ be a compact CW-pair and $\hat{\xi}^k\in \CStr(M\cup H^1,\xi_{H^1})$, $k \in K$, be a family of contact structures such that $\hat{\xi}^k\in \mathcal{G}$ for every $k\in G$. Then, we claim that there is a homotopy of families of contact structures $\hat{\xi}^{k,t}\in \CStr(M\cup H^1,\xi_{H^1})$, $(k,t)\in K\times[0,1]$, such that 
    \begin{itemize}
        \item $\hat{\xi}^{k,t}=\hat{\xi}^k$, for $(k,t)\in K\times\{0\}\cup G\times[0,1]$, and 
        \item $\hat{\xi}^{k,1}\in \mathcal{G}$ for every $k\in K$.
    \end{itemize}
    The claim implies the result. 

    Let us prove the claim. To this end, we will consider the space $X$ of pairs $(\hat{\xi}, e^t)$ where $\hat{\xi}\in\CStr(M\cup H^1,\xi_{H^1})$ and $e^t:D\hookrightarrow M\cup H^1$, $t\in [0,1]$, is a smooth isotopy between the inclusion $i=e^0:D\hookrightarrow M\cup H^1$ of the cocore, and a standard disk embedding  $e^1:D\rightarrow (M\cup H^1,\hat{\xi})$. By Gray stability there is a Serre fibration 
    \[ X\rightarrow \CStr(M\cup H^1,\xi_{H^1}), (\hat{\xi},e^t)\mapsto \hat{\xi}. \]
    Theorem \ref{thm:standardDisks} implies that the fiber of this fibration is weakly contractible.
    
    Therefore, the lift of the map
\[
G\rightarrow \CStr(M\cup H^1,\xi_{H^1}), \qquad
k\mapsto \hat{\xi}^k,
\]
given by the constant isotopy,
\[
G\rightarrow X, \qquad k\mapsto (\hat{\xi}^k,e^{k,t}\equiv i),
\]
extends over the whole parameter space $K$, yielding a family
$(\hat{\xi}^k,e^{k,t})\in X$.

    Apply the smooth isotopy extension theorem to the $K$-family of isotopies $e^{k,t}$ to find diffeomorphisms $\varphi^{k,t}\in \Diff(M\cup H^1)$, $(k,t)\in K\times [0,1]$, such that 
    \begin{itemize}
        \item $\varphi^{k,t}=\Id$ for $(k,t)\in K\times\{0\} \cup G\times [0,1]$, and
        \item $\varphi^{k,t}\circ i= e^{k,t}$.
    \end{itemize}

    The required homotopy of contact structures is then defined as \[ \hat{\xi}^{k,t}:= (\varphi^{k,t})^* \hat{\xi}^k, \qquad (k,t)\in K\times [0,1]. \] This proves the claim and finishes the case of contact $1$-handles. 

   We conclude with contact $3$-handles. Theorem \ref{thm:Standard_spheres} is equivalent to the statement that the space of standard spheres fixed near a point $p$ is weakly homotopy equivalent to the corresponding space of smooth spheres fixed near the same point. Therefore, the same argument as in the case of $1$-handles applies and yields a weak homotopy equivalence \[ \CStr(M\cup H^3,\xi_{H^3},\rel p)
\cong
\CStr(M\cup H^3,\xi_{H^3},\rel S)
=
\CStr(M,\xi)\times \CStr(H^3,\xi_\std). \] The result now follows from Theorem \ref{thm:ContactStructuresR3StandardInfinity}. We leave the details to the reader.
\end{proof}

\subsection{Standard and multistandard tight handlebodies}\label{subsec:StandardHandlebody}

Let $(H_g,\xi)$ be any genus $g$ handlebody equipped with a tight contact structure such that 
\begin{itemize}
    \item the boundary of $H_g$ is convex with dividing set $\Gamma$ and 
    \item it admits a family of pairwise disjoint non-separating disks $\D_1,\ldots,\D_g$ such that $\#\left(\D_i\cap\Gamma\right)=2$.
\end{itemize} 
We say that $(H_g,\xi)$ is a \emph{standard genus $g$ contact handlebody}.

Our gluing theorem, Theorem \ref{thm:Gluing}, allows us to prove our result about standard contact handlebodies stated in the Introduction:

\begin{proof}[Proof of Theorem \ref{thm:StandardHandlebodyINTRO}]
It is a well-known fact, see e.g. \cite{ColinBourbaki}, that, after a slight perturbation along the boundary, a standard tight handlebody can be obtained from a tight ball by attaching several contact $1$-handles. We provide the argument for completeness. Let $(H_g,\xi)$ be a standard contact handlebody. It follows from the Legendrian Realization Principle \cite{Honda} that we may assume that the boundary of each non-separating disk is Legendrian. Proposition \ref{prop:GirouxChangingBoundary} implies that these operations do not change the homotopy type of $\CStr(H_g,\xi)$. Since $(H_g,\xi)$ is tight, the boundary of each disk has negative $\tb$, and thus we can perturb the disks to be convex (Theorem \ref{thm:GirApp}). Since $(H_g,\xi)$ is tight, the dividing set of each disk has no closed components (Theorem \ref{thm:GirouxTightnessCriterion}). In particular, Lemma \ref{lem:DividingSetIntersection} implies that the dividing set must consist of a single segment joining two boundary points. Apply the Realization Theorem again to make the non-separating disks standard. This proves that any tight contact structure on $H_g$ can be homotoped, relative to the boundary, in such a way that the non-separating disks are standard, i.e. the contact handlebody is obtained from a $3$-ball by attaching several contact $1$-handles.

It follows from this observation, Theorem \ref{thm:ContactStructuresR3StandardInfinity}, and Theorem \ref{thm:Gluing} for contact $1$-handles that $\CStr(H_g,\xi)$ is weakly contractible. Since $\Diff(H_g)$ is connected (and weakly contractible) by Corollary \ref{cor:DiffHandlebody} this implies that $\Cont(H_g,\xi)\hookrightarrow \Diff(H_g)$ is a weak homotopy equivalence. This concludes the proof. 
\end{proof}

We have just seen that, after a slight perturbation, a standard contact handlebody admits a family of standard convex non-separating disks. Similarly, we define a \emph{multistandard handlebody} to be a tight handlebody admitting a family of convex non-separating disks with Legendrian boundary and $\partial$-parallel dividing set. Our description of the space of multistandard disks (Theorem \ref{thm:MultiStandardDisks}) allows Theorem \ref{thm:StandardHandlebodyINTRO} to be generalized to multistandard handlebodies. The statement is as follows:

\begin{theorem}\label{thm:MultistandardHandlebody}
    Let $(H_g,\xi)$ be a multistandard genus $g$ handlebody. Then the space $\CStr(H_g,\xi)$ is weakly contractible. Equivalently, the inclusion $\Cont(H_g,\xi)\hookrightarrow \Diff(H_g)$ is a weak homotopy equivalence. In particular, the group $\Cont(H_g,\xi)$ is weakly contractible.
\end{theorem}
\begin{proof}
It is enough to prove the second statement. Let $e:\D^2\hookrightarrow(H_g,\xi)$ be an embedding of a multistandard non-separating disk with characteristic foliation $\mathcal{F}$. Consider the following commutative diagram in which the rows are fibrations:
 \begin{displaymath} 
\xymatrix@M=10pt{
\Cont(H_{g-1},\xi) \ar@{^{(}->}[r]\ar@{^{(}->}[d] &\Cont(H_g,\xi)\ar[r]\ar@{^{(}->}[d] &\Emb ((\D^2,\mathcal{F}) ,(H_g,\xi),e) \ar@{^{(}->}[d] \\
\Diff(H_{g-1}) \ar@{^{(}->}[r] &\Diff(H_g)\ar[r] &\Emb(\D^2,H_g,e).}
\end{displaymath}
Since the right vertical arrow is a weak homotopy equivalence (Theorem \ref{thm:MultiStandardDisks}), it follows that the inclusion $\Cont(H_g,\xi)\hookrightarrow\Diff(H_g)$ is a weak homotopy equivalence if and only if the inclusion $\Cont(H_{g-1},\xi)\hookrightarrow \Diff(H_{g-1})$ is a weak homotopy equivalence. Moreover, the inclusion $\Cont(H_0,\xi)\hookrightarrow \Diff(H_0)$ is a weak homotopy equivalence because of Corollary \ref{coro:ContactomorphismsTightBall}, so the result follows by induction. The weak contractibility statement follows directly from Corollary \ref{cor:DiffHandlebody}.
\end{proof}

A useful consequence of Theorem \ref{thm:StandardHandlebodyINTRO} concerns $(J^1\NS^1,\xi_\std)$, the $1$-jet space over the circle. Recall that there is a contactomorphism between $(J^1\NS^1,\xi_\std)$ and the space $(\NS(T^*\R^2),\xi)$ of cooriented contact elements over $\R^2$, where $\xi=\ker(\cos\theta\, dx-\sin\theta\, dy)$, $(\theta,x,y)\in \NS^1\times\R^2=\NS(T^*\R^2)$. This contactomorphism is explicitly given by
\[
f:(\NS(T^*\R^2),\xi)\longrightarrow (J^1\NS^1,\xi_\std),
\qquad
(\theta,x,y)\longmapsto
(\theta,-x\sin\theta-y\cos\theta,x\cos\theta-y\sin\theta).
\]

\begin{corollary}\label{cor:J1S1}
\begin{itemize}
    \item [(i)] The space of tight contact structures on $\NS^1\times\R^2$ that coincide with $\xi_\std$ at infinity is  weakly contractible. 
    \item [(ii)] The group of compactly supported  contactomorphisms of $(J^1\NS^1,\xi_\std)$ is weakly homotopy equivalent to the corresponding diffeomorphism group and, in particular, is weakly contractible. 
    \item [(iii)] The natural inclusion $\Diff(\D^2)\hookrightarrow\Cont(\NS(T^*\D^2),\xi)$ is a weak homotopy equivalence.
\end{itemize}
\end{corollary}
\begin{proof}
Since the diffeomorphisms $h_t(\theta,y,z)=(\theta,ty,tz)$, $t>0$, are contactomorphisms of $(J^1\NS^1,\xi_\std)$, we may assume that the ambient manifold is a solid torus $\NS^1\times\D^2\subseteq \NS^1\times \R^2=J^1\NS^1$ with the restricted contact structure. This contact structure is clearly a standard genus $1$ handlebody, so (i) and (ii) follow from Theorem \ref{thm:StandardHandlebodyINTRO}. (iii) follows directly from (ii) and the fact that $\Diff(\D^2)$ is weakly contractible because of Smale's Theorem \cite{Smale}.
\end{proof}

\begin{remark} \label{rem:Girouxuniquehandle}
The previous result at the level of path-connected components is a well-known theorem of E. Giroux \cite{Giroux}.
\end{remark}

\subsection{The standard Legendrian unknot in the standard contact $3$-sphere}

We can apply the previous result to compute the weak homotopy type of the space of standard Legendrian unknots in the standard contact $3$-sphere. 

We start with a useful observation about Legendrian embedding spaces in $(\NS^3,\xi_\std)$. Recall that $\xi_\std$ is defined as the complex tangencies of $\NS^3$ with respect to the standard complex structure in $\C^2$. Without loss of generality, we assume that all our Legendrian embeddings have normalized derivative at $t=0$. Then, there is a natural fibration \[ \Leg_{N,jN} (\NS^3,\xi_\std)\hookrightarrow \Leg(\NS^3,\xi_\std)\rightarrow \U(2)\equiv \NS(\xi_\std), \qquad \Lambda\mapsto A_{\Lambda}:= (\Lambda(0),\Lambda'(0) ) \] where \[\Leg_{N,jN}(\NS^3,\xi_\std)=\{\Lambda\in\Leg(\NS^3,\xi_\std): \Lambda(0)=N,\Lambda'(0)=jN\}\] is the space of long Legendrian embeddings in $(\NS^3,\xi_\std)$. This space is weakly homotopy equivalent to the subspace of Legendrian embeddings which coincide with the standard Legendrian great circle $\Lambda_{-1,0}(t)=(\cos t, \sin t)$ near $t=0$. 

Since $\U(2)\subseteq \Cont(\NS^3,\xi_\std)$ acts by contactomorphisms on the sphere, this yields a homotopy equivalence \begin{equation}\label{eq:HomotopyDecompositionLegendrians}
\Phi: \Leg(\NS^3,\xi_\std)\rightarrow \U(2)\times \Leg_{N,jN}(\NS^3,\xi_\std), \qquad \Lambda\mapsto (A_\Lambda, A_{\Lambda}^{-1} \Lambda)
\end{equation} with homotopy inverse \begin{equation}\label{eq:InverseHomotopyDecomposition}  \Psi:\U(2)\times \Leg_{N,jN} (\NS^3,\xi_\std) \rightarrow \Leg(\NS^3,\xi_\std), \qquad (A,\Lambda)\mapsto A\Lambda.
\end{equation}

With these decompositions at hand we are ready to determine the weak homotopy type of $\Leg(\NS^3,\xi_\std,\Lambda_{-1,0})$, the space of standard Legendrian unknots in the standard $3$-sphere. 

\begin{proof}[Proof of Theorem \ref{thm:UnknotsS3Intro}(i)]
    Observe that the space of parametrized Legendrian great circles is canonically identified with $V_{4,2}^{\Lag}\equiv \U(2)$ and, under this identification, the natural inclusion $V_{4,2}^{\Lag}\hookrightarrow \Leg(\NS^3,\xi_\std,\Lambda_{-1,0})$ is precisely the restriction of (\ref{eq:InverseHomotopyDecomposition}) to $\U(2)\times \{\Lambda_{-1,0}\}$. Therefore, it is enough to prove that $\Leg_{N,jN} (\NS^3,\xi_\std, \Lambda_{-1,0})$ is weakly contractible. 

    For this we consider the fibration 
    \[ \Cont(C(\Lambda_{-1,0},\NS^3),\xi_\std)\hookrightarrow \Cont(\NS^3,\xi_\std,\rel N)\cong \Cont(\D^3,\xi_\std)\rightarrow \Leg_{N,jN} (\NS^3,\xi_\std, \Lambda_{-1,0}). \]
 The total space of the fibration is weakly contractible by Corollary \ref{coro:Contcontract}. Moreover, since the complement of $\Lambda_{-1,0}$ is contactomorphic to $(J^1\NS^1,\xi_\std)$, see e.g. \cite{DingGeigesUnknot}, the fiber can be identified with the group of compactly supported contactomorphisms of $(J^1\NS^1,\xi_\std)$ and is therefore also weakly contractible by Corollary \ref{cor:J1S1}. It follows that $\Leg_{N,jN}(\NS^3,\xi_\std,\Lambda_{-1,0})$ is weakly contractible.
\end{proof}

The previous result, together with the contractibility of the space of Lagrangian fillings of the Legendrian unknot $\Lambda_{-1,0}$ in $(\D^4,\lambda_\std)$ proved by Eliashberg--Polterovich \cite{EliashbergPolterovich}, readily implies the following. Denote by $\Lag^0(\D^2,(\D^4,\lambda_\std))$ the space of proper Lagrangian embeddings of the disk $\D^2$ that bound a Legendrian unknot with $\tb=-1$; i.e. the space of exact Lagrangian fillings of \emph{all} Legendrian unknots with $\tb=-1$.

\begin{corollary}\label{cor:Lagrangians}
The natural inclusion $\U(2)\hookrightarrow \Lag^0(\D^2,(\D^4,\lambda_\std))$ is a weak homotopy equivalence. 
\end{corollary}

\subsection{Contact Cell Decompositions and the space of Giroux skeletons}

We recall the contact cell decomposition of a contact $3$-manifold $(M,\xi)$ introduced by E. Giroux \cite{GirouxICM} in his proof of the celebrated correspondence between open books, up to stabilization, and contact structures in dimension $3$.

\begin{definition}[Giroux]\label{def:ContactCell}
	Let $(M,\xi)$ be a compact contact $3$-manifold. A \emph{contact cell decomposition}  of $(M,\xi)$ is a cell decomposition satisfying the following properties:
	\begin{enumerate}
		\item [(i)] The $1$--skeleton is a Legendrian graph,
		\item [(ii)] Each $3$--cell is a Darboux ball, and
		\item [(iii)] Each $2$--cell is attached to a Legendrian unknot with $\tb=-1$.
	\end{enumerate}
\end{definition}

Giroux proved that every contact $3$-manifold admits such a contact cell decomposition, which in turn implies the existence of an adapted open book decomposition \cite{GirouxICM} and of a contact handle decomposition, obtained by thickening the contact cells \cite{Ozbagci}.

Let $(M,\xi)$ be a compact contact $3$-manifold equipped with a contact cell decomposition. Denote by $G_i\subseteq (M,\xi)$ the $i$-skeleton. The $1$-skeleton $G_1$ is an embedded Legendrian graph in $(M,\xi)$, which we call the \emph{Giroux $1$-skeleton}. In this subsection, we show that, in a precise sense, the contact-topological information of $(M,\xi)$ is encoded by $G_1$. Let $GH_1:=\overline{\Op(G_1)}$ be the handlebody given by a standard neighborhood of $G_1$. It follows that $(GH_1,\xi)$ and $(GH_2,\xi):=(\overline{M\backslash GH_1},\xi)$ determine a Heegaard splitting of $(M,\xi)$ in which each handlebody is a standard contact handlebody. See \cite{ColinBourbaki,GirouxICM}. As a consequence of Theorem \ref{thm:StandardHandlebodyINTRO}, we obtain the following lemma.

\begin{lemma}\label{lem:GirouxDecomposition}
Let $(M,\xi)$ be a compact contact $3$-manifold equipped with a contact cell decomposition. Then, the inclusion $\Cont(M,\xi;\rel G_1)\hookrightarrow\Diff(M;\rel G_1)$ is a weak homotopy equivalence. In particular, the group $\Cont(M,\xi;\rel G_1)$ is weakly contractible.
\end{lemma}

Fix a parametrization $g_1:G_1\hookrightarrow (M,\xi)$ of the Giroux $1$-skeleton. Denote by $\Leg^{SK}(M,\xi)$ the space of embeddings of Legendrian graphs that coincide with $g_1$ on $\Op(G_0)$ and are Legendrian isotopic to $g_1$ relative to $\Op(G_0)$. We also denote by $\Emb^{SK}(M)$ the smooth counterpart of $\Leg^{SK}(M,\xi)$. It then follows from the previous lemma and Corollary \ref{cor:DiffHandlebody} that

\begin{corollary}\label{cor:GirouxSkeleton}
Let $(M,\xi)$ be a compact contact $3$-manifold. Then, there is a commutative diagram \begin{displaymath} 
\xymatrix@M=10pt{
\Cont_0(M,\xi;\rel G_0)\ar[r]\ar@{^{(}->}[d] &\Leg^{SK}(M,\xi) \ar@{^{(}->}[d] \\
\Diff_0(M;\rel G_0)\ar[r] &\Emb^{SK}(M),}
\end{displaymath}
where the horizontal arrows are weak homotopy equivalences. 
\end{corollary}

Thus, at the level of contactomorphism groups, the difference between contact and smooth topology is entirely encoded by the space of Legendrian embeddings of the Giroux $1$-skeleton.

\section{The Legendrian unknot in a tight contact $3$-manifold} \label{sec:LegendrianUnknots}

In this section, we generalize Theorem \ref{thm:UnknotsS3Intro}(i) to any tight contact $3$-manifold and to any Legendrian unknot $\Lambda$ satisfying $\tb(\Lambda)+|\rot(\Lambda)|=-1$, thereby proving Theorem \ref{thm:LegendrianAndTransverseUnknotsINTRO}(i). To this end, we establish a weak homotopy equivalence between the based loop space of long Legendrian unknots and the space of convex Seifert disks of a fixed unknot. Finally, we indicate the proofs of Theorems \ref{thm:LegendrianAndTransverseUnknotsINTRO}(ii) and \ref{thm:UnknotsS3Intro}(ii), concerning transverse unknots with self-linking number $\operatorname{sl}=-1$. 

\subsection{Legendrian unknots and Seifert disks}

Let $\Lambda:\NS^1\hookrightarrow (M,\xi)$ be a Legendrian unknot in a contact $3$-manifold $(M,\xi)$ such that $(\Lambda(0),\Lambda'(0))=(p,v)$. Assume that there exists an embedding $e:\D^2\rightarrow(M,\xi)$ of a convex Seifert disk for $\Lambda$ such that $e(1,0)=\Lambda(0)$ and $\Op(e(\D^2))$ is tight. Denote by $\mathcal{F}=e^*\xi$ the characteristic foliation induced by $e$.

We consider the space $\Emb((\D^2,\mathcal{F}),(M,\xi);\frac12\partial)$ of convex disk embeddings with characteristic foliation $\mathcal{F}$ that are Seifert disks for some long Legendrian unknot in the component $\Leg_{(p,v)}(M,\xi,\Lambda)$ and that coincide with $e$ on an open neighborhood of $(1,0)\in\D^2$. This space clearly contains the subspace $\Emb((\D^2,\mathcal{F}),(M,\xi))$ of convex disk embeddings with fixed boundary $\Lambda$ and characteristic foliation $\mathcal{F}$.

These spaces fit into a fibration
\[
\Emb((\D^2,\mathcal{F}),(M,\xi))
\hookrightarrow
\Emb((\D^2,\mathcal{F}),(M,\xi);\frac12\partial)
\longrightarrow
\Leg_{(p,v)}(M,\xi,\Lambda).
\]
Our main observation is that the space of (tight) convex Seifert disks of long Legendrian unknots with half of the boundary fixed is weakly contractible. This is the content of the following result.

\begin{lemma}\label{lem:UnknotsVSDisks}
The space $\Emb((\D^2,\mathcal{F}),(M,\xi);\frac12\partial)$ is weakly contractible. In particular, there is a weak homotopy equivalence \[
\Omega\Leg_{(p,v)}(M,\xi,\Lambda)
\cong
\Emb((\D^2,\mathcal{F}),(M,\xi)).
\]
\end{lemma}
\begin{proof}
Let $\hat{e}:\NS^2\rightarrow(M,\xi)$ be the convex sphere obtained by gluing (and rounding the corners) a copy of $e(\D^2)$ with itself in a vertically invariant neighborhood of the disk. Denote by $\mathcal{F}_2=\hat{e}^*\xi$ the characteristic foliation of $\hat{e}$ and consider the space $\Emb((\NS^2,\mathcal{F}_2),(M,\xi), \rel p)$ of embeddings of convex spheres with characteristic foliation $\mathcal{F}_2$ into $(M,\xi)$ that are fixed near $p=e(1,0)$ and also bound a tight $3$-ball. We claim that $\Emb((\NS^2,\mathcal{F}_2),(M,\xi), \rel p)$ is weakly contractible. To see this consider the space $\CEmb((\D^3,\xi),(M,\xi), \rel p)$ of contact embeddings of tight balls with boundary characteristic foliation $\mathcal{F}_2$ and fixed near $p$. There is a fibration \[ \Cont(\D^3,\xi)\hookrightarrow \CEmb((\D^3,\xi),(M,\xi), \rel p)\rightarrow \Emb((\NS^2,\mathcal{F}_2),(M,\xi), \rel p) \]
The fiber is weakly contractible because of Corollary \ref{coro:ContactomorphismsTightBall}. On the other hand, the total space is weakly homotopy equivalent to the space of Darboux balls fixed near $p$ and therefore also weakly contractible by the contact Alexander trick. Therefore, the base $\Emb((\NS^2,\mathcal{F}_2),(M,\xi), \rel p)$ is weakly contractible. 

There is a second fibration 
\[
\Emb((\NS^2,\mathcal{F}_2),(M,\xi),\rel \D^2)
\hookrightarrow
\Emb((\NS^2,\mathcal{F}_2),(M,\xi),\rel p)
\longrightarrow
\Emb((\D^2,\mathcal{F}),(M,\xi);\frac12\partial).
\]
The fiber consists of those embeddings fixed on one hemisphere. We claim that the fiber $\Emb((\NS^2,\mathcal{F}_2),(M,\xi),\rel \D^2)$ is also weakly contractible, from which the result follows. Indeed, using a fibration analogous to the previous one, one sees that every compact family of spheres in
$\Emb((\NS^2,\mathcal{F}_2),(M,\xi),\rel\D^2)$
extends to a compact family of tight ball embeddings. We may assume that these balls are vertically invariant neighborhoods of the fixed hemisphere. This allows us to contact isotope all the non-fixed hemispheres into a common ball, arbitrarily close to the fixed hemisphere, and then apply Lemma \ref{lem:DisksBall}.
\end{proof}

\begin{remark}
    Observe that we do not require $(M,\xi)$ to be tight in this result.
\end{remark}

\subsection{Proof of Theorem \ref{thm:LegendrianAndTransverseUnknotsINTRO}(i)}

Consider a long Legendrian unknot $\Lambda$ with $\tb(\Lambda)+|\rot(\Lambda)|=-1$ in a tight contact manifold $(M,\xi)$. 

Consider the commuting diagram in which the rows are fibrations \begin{displaymath} 
\xymatrix@M=10pt{
\Emb((\D^2,\mathcal{F}),(M,\xi)) \ar[d] \ar[r]& \Emb((\D^2,\mathcal{F}),(M,\xi); \frac{1}{2}\partial) \ar[d] \ar[r] & \Leg_{(p,v)}(M,\xi,\Lambda)  \ar[d] \\
\Emb(\D^2,M)  \ar[r] & \Emb(\D^2,M,\frac{1}{2}\partial)   \ar[r] & \Emb_{(p,v)}(\NS^1,M,\Lambda).
}
\end{displaymath}
Here, $\Emb(\D^2,M,\frac{1}{2}\partial)$ is the smooth counterpart of the space $\Emb((\D^2,\mathcal{F}),(M,\xi); \frac{1}{2}\partial)$. Observe that this space is weakly contractible. Moreover, $\Emb((\D^2,\mathcal{F}),(M,\xi); \frac{1}{2}\partial)$ is weakly contractible by Lemma \ref{lem:UnknotsVSDisks}, and hence the middle vertical arrow is a weak homotopy equivalence.

Finally, the condition $\tb(\Lambda)+|\rot(\Lambda)|=-1$  implies that $(\D^2,\mathcal{F})$ is a multistandard disk. Therefore, the left vertical arrow is also a weak homotopy equivalence by Theorem \ref{thm:MultiStandardDisks}. The result follows. \qed

\subsection{The transverse unknot with self-linking number $\operatorname{sl}=-1$: Proof of Theorems \ref{thm:LegendrianAndTransverseUnknotsINTRO}(ii) and \ref{thm:UnknotsS3Intro}(ii)}

In this final subsection, we indicate how to prove Theorems \ref{thm:LegendrianAndTransverseUnknotsINTRO}(ii) and \ref{thm:UnknotsS3Intro}(ii). Since the argument is completely analogous to the Legendrian case, we only outline the main steps and leave the details to the reader.

To prove Theorem \ref{thm:LegendrianAndTransverseUnknotsINTRO}(ii), one can establish the analogue of Lemma \ref{lem:UnknotsVSDisks} for disks with transverse boundary. The result then follows from our theorem on mini-disks (Theorem \ref{thm:MiniDisks}).

To deduce Theorem \ref{thm:UnknotsS3Intro}(ii), consider the fibration
\[
\Trans_{N,iN}(\NS^3,\xi_\std,T_{-1})
\hookrightarrow
\Trans(\NS^3,\xi_\std,T_{-1})
\longrightarrow
\NS^3=\SU(2).
\]
There is a canonical identification of the base with the space $V_{4,2}^{\C}\equiv\SU(2)$ of parametrized transverse great circles. Therefore, the result follows by combining Theorems \ref{thm:LegendrianAndTransverseUnknotsINTRO}(ii) and \ref{thm:SmoothLongUnknots}.
\qed


\section{The standard tight $\NS^1\times \NS^2$}\label{sec:S1xS2}


In this section, we prove Theorem \ref{thm:ContactomorphismsS1xS2} for the standard tight contact manifold $(\NS^1\times\NS^2,\xi_\std)$. The contact mapping class group was first studied in detail by Ding and Geiges in \cite{DingGeiges}, and was later completely determined, independently of our work, by Min in \cite{HyunkiMin}. The main ingredient in our proof is Theorem \ref{thm:Standard_spheres}.

\subsection{Standard non-separating spheres in $\NS^1\times \NS^2$}

We begin with a lemma.
\begin{lemma}\label{lem:NonSeparatingS2}
    Let $e:\NS^2\hookrightarrow\{0\}\times\NS^2\subseteq (\NS^1\times\NS^2,\xi_\std)$ be the standard convex inclusion. Then, the space $\Emb((\NS^2,\mathcal{F}_\std),(\NS^1\times\NS^2,\xi_\std),e)$ is weakly homotopy equivalent to $\NS^1\times \U(1)$. 
\end{lemma}
\begin{proof}
The space $\Emb(\NS^2,\NS^1\times\NS^2,e)$ is weakly homotopy equivalent to $\NS^1\times\SO(3)$. This follows from Hatcher's computation of the homotopy type of the diffeomorphism group of $\NS^1\times\NS^2$ (see \cite{HatcherSpheres}). Consider the following commutative diagram 
	\begin{displaymath} 
\xymatrix@M=10pt{
	\NS^1\times\U(1) \ar@{^{(}->}[r]\ar@{^{(}->}[d]& \NS^1\times\SO(3) \ar@{^{(}->}[d] \\
	\Emb((\NS^2,\mathcal{F}_\std),(\NS^1\times\NS^2,\xi_\std),e)\ar@{^{(}->}[r]& \Emb(\NS^2,\NS^1\times\NS^2,e). }
\end{displaymath}
The statement now follows easily from Theorem \ref{thm:Standard_spheres} and the Five Lemma.
\end{proof}

\subsection{Proof of Theorem \ref{thm:ContactomorphismsS1xS2}}

We first prove (ii). Let $\mathcal{S}_\std$ be the space of embeddings of non-separating standard spheres into $(\NS^1\times\NS^2,\xi_\std)$, and let $\mathcal{S}$ denote its smooth counterpart. The space $\mathcal{S}$ has two path-connected components: the one containing the inclusion $e:\NS^2\hookrightarrow \{0\}\times\NS^2\subseteq \NS^1\times\NS^2,$
and the one containing $e\circ A$, where $A:\NS^2\rightarrow\NS^2$ denotes the antipodal map.

Since a standard convex sphere has an antipodal symmetry, the embedding $e\circ A$ is also standard. Moreover, it follows from the work of Colin \cite{Colin} that the inclusion $\mathcal{S}_\std\hookrightarrow \mathcal{S}$ induces an isomorphism on $\pi_0$. Thus,
\[
\mathcal{S}_\std=
\Emb((\NS^2,\mathcal{F}_\std),(\NS^1\times\NS^2,\xi_\std),e)
\cup
\Emb((\NS^2,\mathcal{F}_\std),(\NS^1\times\NS^2,\xi_\std),e\circ A).
\]
Therefore, Lemma \ref{lem:NonSeparatingS2} implies that the natural inclusion
\[
\NS^1\times\Ort(2)\hookrightarrow\mathcal{S}_\std
\]
is a weak homotopy equivalence.

We prove (ii) by studying the fibration
\[
\Cont(I\times \NS^2,\xi_\std)
\hookrightarrow
\Cont(\NS^1\times\NS^2,\xi_\std)
\rightarrow
\mathcal{S}_\std,
\]
where the projection is given by postcomposition of the embedding $e$ with a contactomorphism. Notice that this projection is surjective. Indeed, the diffeomorphism
\[
\varphi_A:\NS^1\times\NS^2\rightarrow\NS^1\times\NS^2,
\qquad
(\theta,p)\mapsto(-\theta,A(p)),
\]
is a contactomorphism of $(\NS^1\times\NS^2,\xi_\std)$.

The fiber of this fibration was determined in Lemma \ref{lem:localSEmb}: the inclusion
\[
\Omega\U(1)\hookrightarrow\Cont(I\times \NS^2,\xi_\std)
\]
is a weak homotopy equivalence.

To conclude, consider the following commuting diagram, whose rows are fibrations:
\begin{displaymath}
\xymatrix@M=10pt{
\Omega \U(1) \ar@{^{(}->}[d]\ar@{^{(}->}[r] &
    \NS^1\times\Ort(2)\times\Omega\U(1) \ar@{^{(}->}[d] \ar[r]&
\NS^1\times\Ort(2) \ar@{^{(}->}[d] \\
\Cont(I\times \NS^2,\xi_\std) \ar@{^{(}->}[r] &
    \Cont(\NS^1\times\NS^2,\xi_\std) \ar[r] &
\mathcal{S}_\std.
}
\end{displaymath}
We have just shown that the maps between the bases and the fibers are weak homotopy equivalences. It follows that the map between the total spaces is also a weak homotopy equivalence. This proves (ii).

We now turn to the proof of (i). We deduce the statement from (ii) by studying the Gray fibration
\[
\Cont_0(\NS^1\times\NS^2,\xi_\std)
\rightarrow
\Diff_0(\NS^1\times\NS^2)
\rightarrow
\CStr(\NS^1\times\NS^2,\xi_\std).
\]
The group $\Diff_0(\NS^1\times\NS^2)$ is weakly homotopy equivalent to $\NS^1\times L^0\SO(3)$ by \cite{HatcherSpheres}, where $L^0\SO(3)$ denotes the space of contractible loops in $\SO(3)$. On the other hand, $\Cont_0(\NS^1\times\NS^2,\xi_\std)$ is weakly homotopy equivalent to
$ \NS^1\times L^{\even}\U(1)$ by (ii). Here $L^{\even}\U(1)$ denotes the subspace of loops in $\U(1)$ that are contractible when regarded as loops in $\SO(3)$.

To conclude, we describe these weak homotopy equivalences in terms of suitable evaluation maps.

Let $\pi:\NS^1\times\NS^2\longrightarrow\NS^1$ denote the projection, and fix the base point $p_0=(0,(1,0,0))\in\NS^1\times\NS^2$. Recall also that we have fixed a trivialization of the tangent bundle of $\NS^1\times\NS^2$. With respect to this trivialization, the weak equivalence for contactomorphisms is given by
\[
\ev_{\NS^1}:
\Cont_0(\NS^1\times\NS^2,\xi_\std)
\longrightarrow
\NS^1\times L^{\even}\U(1),
\qquad
\varphi\longmapsto
\left(
\pi\circ\varphi(p_0),
(\varphi_*)_{|\NS^1\times\{(1,0,0)\}}
\langle\partial_y,\partial_z\rangle
\right).
\]
Similarly, for diffeomorphisms it is given by
\[
\ev_{\NS^1}:
\Diff_0(\NS^1\times\NS^2)
\longrightarrow
\NS^1\times L^0\SO(3),
\qquad
\varphi\longmapsto
\left(
\pi\circ\varphi(p_0),
(\varphi_*)_{|\NS^1\times\{(1,0,0)\}}
\langle\partial_\theta,\partial_y,\partial_z\rangle
\right).
\]

Thus, the Gray fibration fits into the commuting diagram
\begin{displaymath}
\xymatrix@M=10pt{
\Cont_0(\NS^1\times\NS^2,\xi_\std) \ar^{\ev_{\NS^1}}[d]\ar@{^{(}->}[r] &
    \Diff_0(\NS^1\times\NS^2)\ar^{\ev_{\NS^1}}[d] \ar[r]&
\CStr(\NS^1\times\NS^2,\xi_\std) \ar^{\ev_{\NS^1}}[d] \\
\NS^1\times L^{\even}\U(1) \ar@{^{(}->}[r] &
    \NS^1\times L^0\SO(3) \ar[r] &
L\NS^2,
}
\end{displaymath}
where the bottom row is also a fibration, the left and middle vertical maps are the weak homotopy equivalences described above, and the right vertical map is the evaluation map (\ref{eq:EvaluationMapS1xS2}). The result follows. \qed


\section{Legendrian circle bundles over surfaces with boundary}\label{sec:LegendrianFibrations}


In this section, we prove Theorem \ref{thm:LegendrianBoundary} for a general Legendrian circle bundle $(V,\xi)$ over a compact surface $S$ with non-empty boundary.

We use our microfibration method (Theorem \ref{thm:naked}) to adapt the proof given by Giroux and Massot \cite{GirouxMassot} in the $\pi_0$-case. Their proof of the $\pi_0$-statement is based on reducing the complexity of the $3$-manifold $V$ by cutting along an annulus fibered over a properly embedded non-separating arc, and then applying an inductive argument. The base case of the induction is $S=\D^2$, which was already established in Corollary \ref{cor:J1S1}.

\subsection{The space of fibered annuli}

Let $e:A\hookrightarrow (V,\xi)$ be a convex embedding of an annulus fibered over a properly embedded non-separating arc in $S$. The characteristic foliation of $A$, denoted by $\mathcal{F}$, consists of the Legendrian fibers (rulings) and an even number of transverse Legendrian divides, parallel to the arc over which the annulus is fibered. Between each pair of consecutive Legendrian divides there is a dividing curve.

Denote by $\Emb(A,V,e)$ the space of smooth embeddings that coincide with $e$ near the boundary and are isotopic to $e$, and by $\Emb((A,\mathcal{F}),(V,\xi),e)\subseteq\Emb(A,V,e)$ the subspace of convex embeddings with characteristic foliation $\mathcal{F}$.

A key step in the argument of \cite{GirouxMassot} is the proof that
\[
\pi_1(\Emb(A,V,e),\Emb((A,\mathcal{F}),(V,\xi),e))=0.
\]
We need a multiparametric version of this statement, which is the content of the following lemma.

\begin{lemma}\label{lem:Annuli}
The inclusion
\[
\Emb((A,\mathcal{F}),(V,\xi),e)\hookrightarrow \Emb(A,V,e)
\]
is a weak homotopy equivalence.
\end{lemma}
\begin{proof}
We apply the microfibration method (Theorem \ref{thm:naked}) to $E_0=\Emb((A,\mathcal{F}),(V,\xi),e)$ and $E_1=\Emb(A,V,e)$. The density condition (i) in Theorem \ref{thm:naked} was proved by Giroux and Massot in Proposition 2.4 of \cite{GirouxMassot}, and follows from the fact that a Legendrian fiber cannot be destabilized.

It remains to verify condition (ii) in Theorem \ref{thm:naked}, namely, to prove the statement locally in a neighborhood of $e(A)$.

Let $U$ be a small $I$-invariant neighborhood of $A$. Then $(U,\xi)$ is a standard solid torus. Consider the following commutative diagram:
\begin{displaymath}
\xymatrix@M=10pt{
\Cont(U,\xi)^2 \ar@{^{(}->}[r]\ar@{^{(}->}[d] &
\Cont(U,\xi)\ar[r]\ar@{^{(}->}[d] &
\Emb((A,\mathcal{F}),(U,\xi),e) \ar@{^{(}->}[d] \\
\Diff(U)^2 \ar@{^{(}->}[r] &
\Diff(U)\ar[r] &
\Emb(A, U, e).
}
\end{displaymath}
The inclusion maps between the fibers and between the total spaces are weak homotopy equivalences by Theorem \ref{thm:StandardHandlebodyINTRO}. Therefore, the map between the bases is also a weak homotopy equivalence. This verifies condition (ii), and the result follows from Theorem \ref{thm:naked}.
\end{proof}

\subsection{Proof of Theorem \ref{thm:LegendrianBoundary}}

Since the components of $\Diff(S)$ and $\Diff(V)$ are weakly contractible (see \cite{EarlSchatz,GramainDiff,Ivanov,HatcherHomeo,HatcherSurfaces}), it is enough to prove that the inclusion
\[
\Cont(V,\xi)\hookrightarrow \Diff(V)
\]
induces an isomorphism on $\pi_k$ for every $k>0$, and then invoke \cite{GirouxMassot} for the $\pi_0$-case.

Following \cite{GirouxMassot}, we proceed by induction on
\[
n(S)=-2\chi(S)-\beta(S)=\beta(S)+4g(S)-4,
\]
where $g(S)$ denotes the genus of $S$ and $\beta(S)$ its number of boundary components.

The base case $n(S)=-3$ corresponds to Corollary \ref{cor:J1S1}. Assume now that $n(S)>-3$. Let $e:A\hookrightarrow (V,\xi)$ be a convex embedding of an annulus fibered over a properly embedded non-separating arc in $S$. Consider the following commuting diagram, whose rows are fibrations:
\begin{displaymath}
\xymatrix@M=10pt{
\Cont_0(V',\xi) \ar@{^{(}->}[r]\ar@{^{(}->}[d] &
\Cont_0(V,\xi)\ar[r]\ar@{^{(}->}[d] &
\Emb((A,\mathcal{F}),(V,\xi),e) \ar@{^{(}->}[d] \\
\Diff_0(V') \ar@{^{(}->}[r] &
\Diff_0(V)\ar[r] &
\Emb(A, V, e).
}
\end{displaymath}
Here, $(V',\xi)$ is a Legendrian circle bundle over a surface $S'$ satisfying
\[
n(S')<n(S).
\]

By Lemma \ref{lem:Annuli}, the map between the bases is a weak homotopy equivalence. It follows that there are natural isomorphisms
\[
\pi_k(\Diff_0(V),\Cont_0(V,\xi))
\cong
\pi_k(\Diff_0(V'),\Cont_0(V',\xi)),
\qquad k>0.
\]
The result now follows from the induction hypothesis. \qed


\section{The space of Legendrian parametrized $(n,n)$-torus links with max-$\tb$}\label{sec:TorusLinks} 


In this section, we prove Theorem \ref{thm:LegendrianNNlinks}.

\begin{proof}[Proof of Theorem \ref{thm:LegendrianNNlinks}]
The main geometric observation is that every Legendrian $(n,n)$-torus link with maximal $\tb$-number \[\Lambda=\cup_{i=1}^{n} \Lambda_i: \sqcup_{i=1}^n \NS^1\hookrightarrow (\NS^3,\xi_\std)\] is Legendrian isotopic to $n$ distinct parametrized Legendrian fibers of the Legendrian Hopf fibration $(\NS^3,\xi_\std)\rightarrow \C \PP^1$ induced by the $j$-complex structure on $\R^4$. Consequently, the link exterior $(C(\Lambda,\NS^3),\xi_\std)$ is a Legendrian circle bundle over $S_n$, the $2$-sphere with $n$ holes.

As in (\ref{eq:HomotopyDecompositionLegendrians}), there is a weak homotopy equivalence
\[ \mathcal{L}_n\cong\U(2)\times\mathcal{L}'_n, \]
where $\mathcal{L}'_n\subseteq \mathcal{L}_n$ is the subspace of Legendrian links $\Lambda=\sqcup_{i=1}^n \Lambda_i$ satisfying $\Lambda_1(0)=N$ and $\Lambda'_1(0)=jN$. Thus, it is enough to prove that $\mathcal{L}'_n$ has the weak homotopy type of a $K(\mathcal{M}_n,1)$.

Consider the fibration $$ \Cont(\D^3,\xi_\std)\cong \Cont(\NS^3,\xi_\std,\rel N)\rightarrow \mathcal{L}'_n $$
defined by postcomposition with a fixed Legendrian link. The fiber over a Legendrian link $\Lambda$ can be identified with $\Cont(C(\Lambda,\NS^3),\xi_\std)$.
Since $\Cont(\D^3,\xi_\std)$ is weakly contractible by Corollary \ref{coro:Contcontract}, we obtain a weak homotopy equivalence
\[ \Cont(C(\Lambda,\NS^3),\xi_\std)\cong\Omega \mathcal{L}'_n. \]

By Theorem \ref{thm:LegendrianBoundary}, the natural inclusion \[\Diff(S_{n})\hookrightarrow\Cont(C(\Lambda,\NS^3),\xi_\std)\] is a weak homotopy equivalence. Thus, \[\Omega\mathcal{L}'_n\cong\Diff(S_n).\] Since $\Diff(S_n)$ is homotopically discrete and $\mathcal{M}_n=\pi_0(\Diff(S_n))$ by definition, it follows that $\mathcal{L}'_n$ has the weak homotopy type of a $K(\mathcal{M}_n,1)$. This concludes the argument.
\end{proof}


\section{Parametric families of tight charts}\label{sec:TightCharts}


A classical theorem of Eliashberg \cite{Eliashberg} states that for every tight contact structure $\xi$ on $\R^3$ there exists a contactomorphism $\varphi:(\R^3,\xi)\rightarrow (\R^3,\xi_\std)$. The goal of this section is to generalize this result to the parametric setting:
\begin{theorem}\label{thm:EliashbergParametrico}
	Let $\xi^k$, $k\in K$, be a compact family of tight contact structures on $\R^3$. Assume that there is a choice of contact form $\alpha_k$ on $\Op(0)$ for which the Reeb field satisfies $R_k(0)= \partial_z$. Then, there exists a family of contactomorphisms $\varphi_k:(\R^3,\xi^k)\rightarrow(\R^3,\xi_\std)$, $k\in K$.
\end{theorem}
\begin{proof}
The condition at the origin implies that $\xi^k$ may be chosen to be fixed on $\Op(0)$. Since the space of contact structures under consideration retracts, by Alexander's trick, onto the standard contact structure, we may assume that we have a cone of tight contact structures $\xi^{k,t}$, $(k,t) \in K \times [0,1]$, such that $\xi^{k,1}= \xi^k$ and $\xi^{k,0}=\xi_\std$.

Consider the exhaustion $\R^3=\cup_n \D^3(n)$. Observe that $\NS^2(n)=\partial\D^3(n)$ is convex for the standard contact structure $(\R^3,\xi_\std)$. Denote the natural inclusion $\NS^2(n)\hookrightarrow \R^3$ by $e_{0,n}$. 

The condition on the Reeb field and  Theorem \ref{thm:Standard_spheres} imply the existence of a family of embeddings $e_{k,t,n}:\NS^2\rightarrow \R^3 $, $(k,t,n)\in K \times [0,1] \times \Z_+$, such that 
\begin{itemize}
    \item [(i)] $e_{k,t,n}=e_{0,n}$ for $(k,t,n)\in K\times \{0\}\times \Z_+$,
    \item [(ii)] $e_{k,t,n}$ is standard for $\xi^{k,t}$ for every $(k,t,n)$,
    \item [(iii)] The image of $e_{k,t,n}$ is arbitrarily close to $\NS^2(n)$ for $(k,t,n)\in K\times [0,1]\times \Z_+$.
\end{itemize}

For every $n\in \Z_+$, apply the smooth isotopy extension theorem to find a family of diffeomorphisms $F_{k,t,n}:\R^3\rightarrow \R^3$, $(k,t)\in K\times [0,1]$, such that $e_{k,t,n}=F_{k,t,n}\circ e_{0,n}$. By (i) we may assume that $F_{k,t,n}=\Id$ for every $(k,t,n)\in K \times\{0\}\times\Z_+$. Property (iii) implies that we may further assume that for every fixed $n\in \Z_+$ the support of $F_{k,t,n}$ lies inside $\Op(\NS^2(n))$. In particular, for $n\neq m$, the supports of $F_{k,t,n}$ and $F_{k,t,m}$ are disjoint. 

Thus, $F_{k,t} = \cdots \circ F_{k,t,2} \circ F_{k,t,1}$, $(k,t)\in K\times [0,1]$, is a diffeomorphism such that $F_{k,t} \circ e_{0,n}$ is a standard embedding for $\xi^{k,t}$ and, moreover, it is a contactomorphism on small open neighborhoods of the spheres, i.e. $F_{k,t}^* \xi^{k,t}= \xi_{\std}$ on the domain $\bigcup_{n \in \Z^+}\Op(\NS^2(n))$. 

Since $F_{k,t}=\Id$ for every  $(k,t)\in K\times \{0\}$ by construction, we have a homotopy of $K$-families of contact structures $F_{k,t}^*\xi^{k,t}$, $t\in [0,1]$, between $F_{k,1}^*\xi^{k,1}$ and the constant family $F_{k,0}^*\xi^{k,0}=\Id^* \xi_\std=\xi_\std$. Moreover, all these contact structures coincide with $\xi_\std$ on the domain $\bigcup_{n \in \Z^+}\Op(\NS^2(n))$. Therefore,  Gray's flow is complete, since it is confined to the complementary annuli. Thus, there exists a family $G_{k,t}:\R^3\rightarrow \R^3$, $(k,t)\in K\times [0,1]$, such that $(G_{k,t})_*\xi_\std =F_{k,t}^*\xi^{k,t}$. The required family of contactomorphisms is \[ \varphi_k:= (F_{k,1}\circ G_{k,1})^{-1}. \]
\end{proof}

\bibliographystyle{alpha}
\bibliography{biblio}

\end{document}